\documentclass[12pt]{report}

\usepackage{fancyheadings}
\usepackage{graphicx} 
\usepackage{Macros}
\usepackage{tikz-cd}

\DeclareMathOperator\Aut{Aut}
\DeclareMathOperator\proj{proj}
\DeclareMathOperator\GL{GL}
\DeclareMathOperator\Stab{Stab}
\DeclareMathOperator\id{id}

\renewcommand{\theequation}{\arabic{chapter}.\arabic{section}.\arabic{equation}}

\begin{document}

\pagenumbering{roman}
 \addtolength{\headwidth}{\marginparsep}
 \addtolength{\headwidth}{\marginparwidth}

\begin{titlepage}
%
%
%
%
\begin{LARGE}
\begin{bf}
\begin{center}
Holomorphic Flexibility Properties\\ of Spaces of Elliptic Functions
\end{center}
\end{bf}
\end{LARGE}
%
%
\vspace{5.0mm}
%
%
\begin{large}
\begin{center}
David Bowman
\end{center}
%
%
\vspace{2.0mm}
%
%
\begin{em}
\begin{center}
Thesis
submitted for the degree
of\\
Doctor of Philosophy\\
in\\
Pure Mathematics\\
at\\
The University of Adelaide\\
Faculty of Engineering, Computer\\ and Mathematical Sciences\\
\end{center}
\end{em}
%
%
\vspace{4.0mm}
%
%
\begin{center}
School of Mathematical Sciences
\end{center}
%
%
\vspace{4.0mm}
\begin{center} 
\includegraphics[height=2cm]{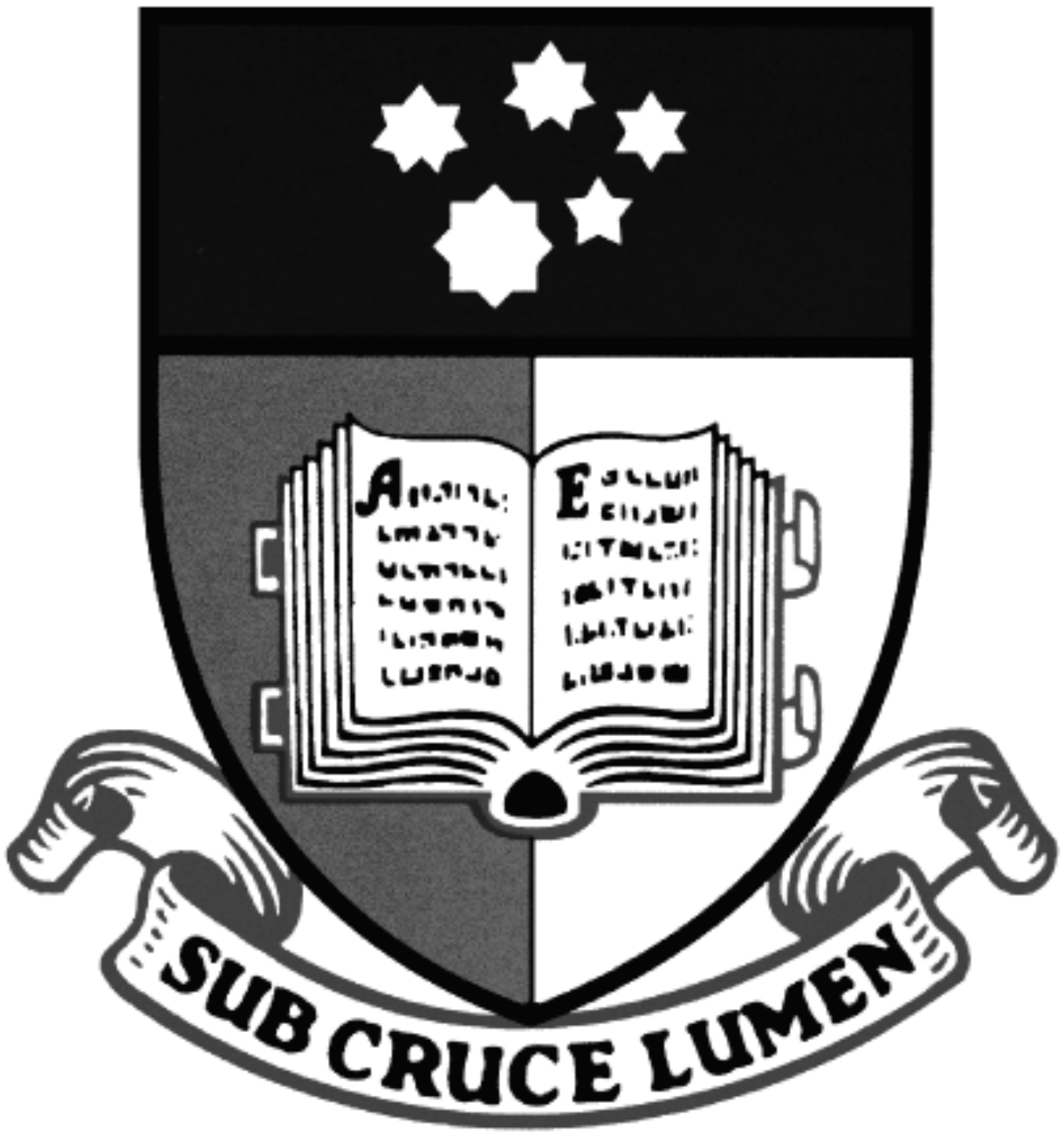}
\end{center}
%
%
%
%
\begin{center}
\today
\end{center}
\end{large}
%
%
\end{titlepage}

\tableofcontents

\chapter*{Signed Statement}
\addcontentsline{toc}{chapter}{Signed Statement}

This work contains no material which has been accepted
for the award of any other degree or diploma in any university or other
tertiary institution and, to the best of my knowledge and belief,
contains no material previously
published or written by another person, except where due reference has
been made in
the text.

\vspace{5mm}

\noindent
I consent to this copy of my thesis, when deposited in the University
Library, being available for loan and photocopying.

\vspace{20mm}

\noindent
{\sc SIGNED}: {\tt .......................} {\sc DATE}: {\tt .......................}

\chapter*{Acknowledgements}
\typeout{Acknowledgements}
\label{ch:acknowledgements}
\addcontentsline{toc}{chapter}{Acknowledgements}

I would like to thank the following for their contributions to this thesis.

Finnur L\'arusson, my principal supervisor, for his unending support, for nurturing the growth of this thesis and for many, many helpful discussions.

Nicholas Buchdahl, my secondary supervisor, for exposing me to many exciting ideas in geometry and analysis.

Franc Forstneri\v c, Philipp Naumann, Richard L\"ark\"ang and Tuyen Truong for many interesting mathematical discussions.

The staff at Kathleen Lumley College for making my stay in Adelaide an enjoyable one.

And finally, my friends and family for their love and support.

\chapter*{Dedication}
\typeout{Dedication}
\label{ch:dedication}
\addcontentsline{toc}{chapter}{Dedication}

To Grandad, for teaching me to count cars.
\chapter*{Abstract}
\typeout{Abstract}
\label{ch:abstract}
\addcontentsline{toc}{chapter}{Abstract}

Let $X$ be an elliptic curve and $\mathbb{P}$ the Riemann sphere. Since $X$ is compact, it is a deep theorem of Douady that the set $\mathcal{O}(X,\mathbb{P})$ consisting of holomorphic maps $X\to \mathbb{P}$ admits a complex structure. If $R_n$ denotes the set of maps of degree $n$, then Namba has shown for $n\geq2$ that $R_n$ is a $2n$-dimensional complex manifold. We study holomorphic flexibility properties of the spaces $R_2$ and $R_3$. Firstly, we show that $R_2$ is homogeneous and hence an Oka manifold. Secondly, we present our main theorem, that there is a $6$-sheeted branched covering space of $R_3$ that is an Oka manifold. It follows that $R_3$ is $\mathbb{C}$-connected and dominable. We show that $R_3$ is Oka if and only if $\mathbb{P}_2\backslash C$ is Oka, where $C$ is a cubic curve that is the image of a certain embedding of $X$ into $\mathbb{P}_2$.

We investigate the strong dominability of $R_3$ and show that if $X$ is not biholomorphic to $\mathbb{C}/\Gamma_0$, where $\Gamma_0$ is the hexagonal lattice, then $R_3$ is strongly dominable.

As a Lie group, $X$ acts freely on $R_3$ by precomposition by translations. We show that $R_3$ is holomorphically convex and that the quotient space $R_3/X$ is a Stein manifold.

We construct an alternative $6$-sheeted Oka branched covering space of $R_3$ and prove that it is isomorphic to our first construction in a natural way. This alternative construction gives us an easier way of interpreting the fibres of the branched covering map.

\pagenumbering{arabic}

\pagestyle{fancy}
\rhead{\thepage}
\cfoot{}
\chapter{Introduction \label{ch:intro}}


\markboth{\uppercase{Chapter}~\protect\ref{ch:intro}. \uppercase{Introduction}}{\uppercase{Chapter}~\protect\ref{ch:intro}. \uppercase{Introduction}}

\section{Background and Context}

Oka theory originates from the principle that certain problems in complex analysis have only topological obstructions. The earliest work in this vein was done by Oka, who showed in 1939 that the second Cousin problem on a domain of holomorphy has a holomorphic solution if and only if it has a continuous solution. It follows that the holomorphic and topological classifications of line bundles over domains of holomorphy are equivalent. Grauert extended this result to $G$-bundles over Stein spaces, where $G$ is a complex Lie group.

The modern development of the theory began with the work of Gromov, who was interested in what is now called the \emph{basic Oka property} (BOP). A complex manifold $X$ is said to satisfy BOP if every continuous map $S\to X$, where $S$ is a Stein manifold, is homotopic to a holomorphic map $S\to X$. Grauert showed that every homogeneous space satisfies BOP. Gromov defined \emph{elliptic manifolds} as a generalisation of homogeneous spaces and proved that not only do they satisfy BOP, but also stronger properties involving deforming continuous families of continuous maps to continuous families of holomorphic maps, as well as approximating such maps.

For basic introductions to Oka theory, see \cite{China} and \cite{Woka?}. For more comprehensive references, see \cite{Forstneric}, \cite{Smack} and \cite{Survey}.

To define an Oka manifold we will use a second approach to Oka theory which arises from the classical theorems of Runge and Weierstrass. Weierstrass' theorem generalises to Cartan's extension theorem, which allows holomorphic functions defined on closed subvarieties of Stein manifolds to be extended to the whole manifold. Runge's theorem, on the other hand, generalises to the Oka-Weil approximation theorem: every holomorphic function on a compact holomorphically convex subset $K$ of a Stein manifold $S$ can be approximated uniformly on $K$ by functions that are holomorphic on all of $S$. These two theorems are classically taken to be properties of Stein manifolds, but they can also been viewed as properties of the target manifold~$\mathbb{C}$.

A complex manifold $X$ is said to satisfy the \emph{interpolation property} (IP) if for every subvariety $S$ of a Stein manifold $Y$, a holomorphic map $S\to X$ has a holomorphic extension to $Y$ if it has a continuous extension.

A complex manifold $X$ is said to satisfy the \emph{approximation property} (AP) if for every holomorphically convex compact subset $K$ of a Stein manifold $Y$, a continuous map $Y\to X$ holomorphic on (a neighbourhood of) $K$ can be uniformly approximated on $K$ by holomorphic maps $Y\to X$.

Cartan's extension theorem and the Oka-Weil approximation theorem then become the statements that $\mathbb{C}$ satisfies IP and AP respectively.

It is not difficult to show that IP implies AP. It is a deep theorem of Forstneri\v c that the converse is also true. Indeed, Forstneri\v c's theorem proves the equivalence of a whole family of related properties. The weakest of these is the \emph{convex approximation property} (CAP). We say that a complex manifold $X$ satisfies CAP if given any convex compact subset $K$ of $\mathbb{C}^n$, every holomorphic map $K\to \mathbb{C}^n$ can be approximated uniformly on $K$ by holomorphic maps $\mathbb{C}^n\to X$.

Any manifold $X$ satisfying one of these equivalent properties is called an \emph{Oka manifold}.

Every elliptic manifold is Oka. Hence every homogeneous manifold is Oka. Furthermore, every Oka manifold satisfies BOP. However, the converse is false. Indeed the unit disc $\mathbb{D}$ is not Oka by Liouville's theorem, but satisfies BOP because it is contractible.

Proving that a manifold is Oka is generally quite difficult. It is necessary to be able to construct new Oka manifolds from old. A crucial result that we will use is a deep theorem that a holomorphic fibre bundle with Oka fibres is Oka if and only if the base is Oka. This result originates from the work of Gromov \cite[Corollary 3.3C$'$]{Gromov}. See also \cite[Theorem 1.2]{Oka Maps}. Note that a special case of this theorem is that an unbranched covering space of a manifold is Oka if and only if the base is Oka.

There are several weaker flexibility properties that all Oka manifolds possess, but are not known to be equivalent to the Oka property. It is often easier to establish these properties than the full Oka property. We say that a complex manifold $X$ is $\mathbb{C}$-\emph{connected} if any two points in $X$ can be joined by a holomorphic map $\mathbb{C}\to X$. We say a manifold of dimension $n$ is \emph{dominable} if for some $x\in X$, there exists a holomorphic map $\mathbb{C}^n\to X$ that is a local biholomorphism at $0$, with $f(0)=x$. If we can find such a map for every $x\in X$, then we say that $X$ is \emph{strongly dominable}.

Now, let $X$ and $Y$ be complex manifolds with $X$ compact. Let $\mathcal{O}(X,Y)$ denote the set of holomorphic maps $X\to Y$ equipped with the compact-open topology. It is a special case of a deep theorem of Douady \cite[Section 10.2]{Douady} that there exists a natural complex structure on $\mathcal{O}(X,Y)$. It is natural in the sense that if $Z$ is any complex space, a map $g:Z\to\mathcal{O}(X,Y)$ is holomorphic if and only if the map $f:Z\times X\to Y$, $f(z,x)=g(z)(x)$, is holomorphic. This gives a natural identification between $\mathcal{O}(Z, \mathcal{O}(X,Y))$ and $\mathcal{O}(Z\times X, Y)$. The complex geometry of such mapping spaces is of intrinsic interest, but is particularly relevant to Oka theory when $Y$ is taken to be Oka. In this case, the study of maps into $\mathcal{O}(X,Y)$ reduces to the study of the space $\mathcal{O}(Z\times X, Y)$, the elements of which can be considered to be holomorphic families of holomorphic maps $Z\to Y$ with compact parameter space $X$. However, the study of Oka properties of such mapping spaces is difficult. In particular, even if $X,Y$ are both manifolds, $\mathcal{O}(X,Y)$ may be singular. See \cite[Section 2.2]{Namba} for an example. This is problematic, as Oka properties do not generalise easily to singular spaces \cite{Rik}. 

The first work in this direction was done by Hanysz who studied the space $\mathcal{O}(\mathbb{P}, \mathbb{P})$ of rational functions in one complex variable \cite{Lex}. Here, $\mathbb{P}$ denotes the Riemann sphere. If $R_n$ denotes the set of holomorphic maps $\mathbb{P}\to\mathbb{P}$ of degree $n$, then we can write
\[
\mathcal{O}(\mathbb{P}, \mathbb{P})= R_0 \sqcup R_1 \sqcup R_2 \sqcup R_3\sqcup \cdots.
\]
The set $R_0$ is the set of constant maps, so is biholomorphic to $\mathbb{P}$, which is well known to be Oka. The set $R_1$ is the M\"obius group, a complex Lie group, so is Oka. Hanysz showed that $R_2$ is a homogeneous space, so is Oka. His main result, using geometric invariant theory, was that $R_3$ is strongly dominable and $\mathbb{C}$-connected.

Our work has been to study $\mathcal{O}(X,\mathbb{P})$, where $X$ is an elliptic curve, that is, a compact Riemann surface of genus $1$.

\section{Results of the Thesis}

Let $X$ be a compact Riemann surface of genus $1$. Then $\mathcal{O}(X,\mathbb{P})$ denotes the set of holomorphic maps $X\to \mathbb{P}$. It is a special case of a theorem of Douady \cite[Section 10.2]{Douady} that $\mathcal{O}(X,\mathbb{P})$ admits a complex structure, uniquely determined by a universal property. Furthermore, if we let $R_n\subset\mathcal{O}(X,\mathbb{P})$ denote the set of degree $n$ holomorphic maps $X\to \mathbb{P}$, then Namba has shown for $n\geq 2$, that $R_n$ is in fact a $2n$-dimensional complex manifold \cite[Proposition 1.1.4]{Namba}.

As in Hanysz's work, $R_0$ is biholomorphic to $\mathbb{P}$, so is Oka. However $R_1=\varnothing$, as any element in $R_1$ would be a biholomorphism.

Our first result is the following, analogous to the result of Hanysz in the genus $0$ case.\\

\e{thm}
The complex manifold $R_2$ is a homogeneous space and hence Oka.
\ee{thm}

Our main theorem is the following.\\

\e{thm}
There exists a $6$-sheeted branched covering space of $R_3$ that is an Oka manifold.
\ee{thm}

It must be noted that it is not known whether the Oka property can be passed down a branched covering map. Hence we are unable to conclude that $R_3$ is Oka. However, our theorem does imply that $R_3$ is $\mathbb{C}$-connected and dominable. Hence our main theorem is similar to that of Hanysz in the degree $3$ case, except that it does not imply strong dominability. 

The proof shall use the following commuting diagram:

\[
\begin{tikzcd}
N\times X\arrow{d}{\lambda\times \text{id}_X}\arrow{rd}{\Psi} \arrow{r}{\zeta}& (N\times X)/S \arrow{d}{\Lambda}&W\arrow{d}{\Phi} \arrow{l}{\tilde{\pi}}\\
\mathbb{P}_2\backslash C\times X \arrow{r}{\tilde{s}_3}\arrow{d}{\proj_2}&R_3/M\arrow{d}{\chi}& R_3\arrow{l}{\pi}\arrow{d}{D_0}\arrow{ld}{\beta}\\
X\arrow{r}{s_3}&X& S^3X \arrow{l}{\psi}
\end{tikzcd}
\]

Let $M$ denote the group of M\"obius transformations. Then $M$ is the automorphism group of $\mathbb{P}$ and acts freely and properly on $R_3$ by postcomposition. The orbit space $R_3/M$ is a complex manifold such that the projection $\pi: R_3\to R_3/M$ is a holomorphic principal $M$-bundle.

Let $S^3 X$ be the $3$-fold symmetric product of $X$. Then $S^3X$ is a $3$-dimensional complex manifold. For $f\in R_3$, let $D_0(f)\in S^3 X$ be the divisor of zeroes of $f$. Then the map $D_0: R_3\to S^3 X$ is holomorphic.

Recall that $X$ may be given the structure of a complex Lie group by identifying it with $\mathbb{C}/\Gamma$, where $\Gamma$ is a lattice in $\mathbb{C}$. Let $\psi: S^3 X \to X$ be the Jacobi map, which takes the divisor $p_1+ p_2 + p_3\in S^3 X$ to the point $p_1+ p_2 + p_3\in X$, where the addition in $S^3 X$ is the formal addition of points, and the addition in $X$ is the Lie group addition. Then $\psi$ is holomorphic, and the composition $\beta=\psi\circ D_0$ is $M$-invariant, so $\beta$ induces a holomorphic map $\chi: R_3/ M \to X$.

Let $s_3: X\to X, z\mapsto 3z$. Then $s_3$ is an unbranched $9$-sheeted covering map. Furthermore, $s_3$ is a normal covering. Namba \cite{Namba} has shown that for a certain smooth cubic curve $C$ in the complex projective plane $\mathbb{P}_2$, biholomorphic to $X$, there exists a map $\tilde{s}_3: \mathbb{P}_2\backslash C\times X\to R_3/ M$ such that the following is a pullback square:

\[
\begin{tikzcd}
\mathbb{P}_2\backslash C \times X \arrow{r}{\tilde{s}_3} \arrow{d}{\text{proj}_2} & R_3/M  \arrow{d}{\chi} \\
X  \arrow{r}{s_3} &  X
\end{tikzcd}
\]
where $\text{proj}_2$ is the projection onto the second component. Then $\tilde{s}_3$ is a $9$-sheeted unbranched covering map $\mathbb{P}_2\backslash C\times X \to R_3/ M$. We might say that $s_3$ unravels the map $\chi$.

Buzzard and Lu \cite[Proposition 5.1]{Blu} have shown that there exists a dominable 6-sheeted branched covering space $N$ of $\mathbb{P}_2\backslash C$. Let $\lambda$ denote the $6$-sheeted branched covering map $N\to \mathbb{P}_2\backslash C$. Hanysz \cite[Proposition 4.10]{Alex} has shown that $N$ is in fact Oka.

Then let $\Psi=\tilde{s}_3\circ(\lambda\times \text{id}_X)$. This map is a $54$-sheeted branched covering map $N\times X \to R_3/M$. Let $S=s_3^{-1}(0)=\Set{x\in X}{3x=0}$. Then $S$ is the group of covering transformations of $s_3$. The action of $S$ on $X$ lifts firstly to a free action on $\mathbb{P}_2\backslash C \times X$ and then to a free action on $N\times X$. The quotient map $\zeta: N\times X\to (N\times X)/S$ is a $9$-sheeted unbranched covering map. Since both $N$ and $X$ are Oka, it follows that $N\times X$ and hence $(N\times X)/S$ are Oka. The map $\Psi$ is $S$-invariant and so induces a $6$-sheeted branched covering map $\Lambda: (N\times X)/S\to R_3/M$.

Taking $W$ to be the pullback of $R_3$ along $\Lambda$, we obtain the following pullback square:
\[
\begin{tikzcd}
W \arrow{r}{\Phi} \arrow{d}{\tilde{\pi}}& R_3 \arrow{d}{\pi}\\
(N\times X)/S \arrow{r}{\Lambda} & R_3/M
\end{tikzcd}
\]
where $\tilde{\pi}: W \to (N\times X)/S$ is the pullback of $\pi$ along $\Lambda$ and $\Phi: W\to R_3$ is the pullback of $\Lambda$ along $\pi$. Then $W$ is smooth since $\pi$ is a submersion and $W$ is a $6$-sheeted branched covering space of $R_3$.

Since $\pi$ is a principal $M$-bundle, it follows that $\tilde{\pi}$ is also a principal $M$-bundle. Since $M$ is a complex Lie group, $M$ is Oka. Hence $W$ is Oka.

Note that it is an open problem whether $\mathbb{P}_2\backslash C$ is Oka \cite[Open Problem B, p.~20]{Survey}. Our work shows that $R_3$ is Oka if and only if $\mathbb{P}_2\backslash C$ is Oka. Indeed, if $R_3$ is Oka, since $\pi$ is a principal bundle with Oka fibre $M$, then $R_3/M$ is Oka. Since $\tilde{s}_3$ is an unbranched covering map, it follows that $\mathbb{P}_2\backslash C\times X$ is Oka. Since a retract of an Oka manifold is easily seen to be Oka, it follows that $\mathbb{P}_2\backslash C$ is Oka. Conversely, if $\mathbb{P}_2\backslash C$ is Oka, then $\mathbb{P}_2\backslash C\times X$ is Oka since the complex Lie group $X$ is Oka. Using the maps $\tilde{s}_3$ and $\pi$, it follows that both $R_3/M$ and $R_3$ are Oka.

Our result implies that $R_3$ is dominable and $\mathbb{C}$-connected. We should note that although this result is analogous to the result of Hanysz, who showed that $R_3$ is strongly dominable in the genus $0$ case, our method of proof is completely different. Hanysz took the categorical quotient of $R_3$ by the action of $M\times M$, acting via pre- and post-composition. In our case, $M$ acts freely on $R_3$ by post-composition, which allows us to take the geometric quotient. This has helped us to obtain our results.

The specific references most crucial to our work are those of Namba \cite{Namba}, Buzzard and Lu \cite{Blu}, and Hanysz \cite{Alex}. Namba's work has been particularly vital. It was Namba who first showed that $R_3$ is a manifold and constructed the map $\tilde{s}_3$. Much of Namba's work is done in the more general case where $X$ is a compact Riemann surface of genus $g\geq 1$. In the genus $g=1$ case, his proofs can be considerably simplified, and we have done this to make our work as self-contained as possible. We have also presented the details of the construction of the map $\lambda: N\to\mathbb{P}_2\backslash C$, which was originally done by Buzzard and Lu and is a crucial ingredient to our proof.

Besides our main theorem, we have proved several additional results.

Let $\Gamma_0$ denote the hexagonal lattice in $\mathbb{C}$. Then by investigating the branching of $\Phi: W\to R_3$ we have shown the following.\\

\e{thm}
If $X$ is not biholomorphic to $\mathbb{C}/\Gamma_0$, then $R_3$ is strongly dominable.
\ee{thm}

The Lie group $X$ acts on $R_3$ by pre-composition. This action is free and so $R_3/X$ is a complex manifold. It is easy to see that $R_3$ is not Stein because the orbits of $X$ are compact subvarieties of $R_3$. However, $R_3$ is holomorphically convex and we have shown the following.\\

\e{thm}
The complex manifold $R_3$ is holomorphically convex and its Remmert reduction is $R_3/X$.
\ee{thm}

Finally, by first taking the quotient by $X$ we have constructed an alternate $6$-sheeted Oka branched covering space $W'$ of $R_3$. We have shown that $W'$ is biholomorphic to $W$ in a natural way. This has allowed us to determine the fibres of the branched covering map $\Phi:W\to R_3$ and provide a geometric interpretation of them.

\section{Further Directions}

Let $X$ be a compact Riemann surface of genus $1$ and $R_3$ the complex manifold of holomorphic maps $X\to \mathbb{P}$ of degree $3$. Let $C$ be a cubic curve in $\mathbb{P}_2$ biholomorphic to $X$.

The most immediate open question arising from our work is to determine whether $\mathbb{P}_2\backslash C$ is strongly dominable in the case when $X$ is biholomorphic to $\mathbb{C}/\Gamma_0$, where $\Gamma_0$ is the hexagonal lattice in $\mathbb{C}$. Our proof of the strong dominability of $\mathbb{P}_2\backslash C$ fails in this case because the map $\lambda: N\to \mathbb{P}_2\backslash C$ has a special kind of branching for this cubic. To resolve the problem, it is necessary to construct a new dominating map $\mathbb{C}^2\to\mathbb{P}_2\backslash C$, but it is not obvious how to do so. If $\mathbb{P}_2\backslash C$ can be shown to be strongly dominable, then the associated $R_3$ would also likely be strongly dominable.

Another direction for further research that follows from our work is to study holomorphic flexibility properties of $R_n$ for $n\geq 4$. Recall that $R_3/M$ has a $9$-sheeted covering space $\mathbb{P}_2\backslash C\times X$. Now, for $n\geq 4$, fix $D\in S^nX$. The Riemann-Roch theorem implies that $\dim H^0(X,\mathcal{O}_D)=n$. Let $\eta_1\dots, \eta_n$ form a basis of $H^0(X,\mathcal{O}_D)$. Then the map $X\to \mathbb{P}_{n-1}$, $x\mapsto [\eta_1(x),\dots, \eta_n(x)]$, is an embedding of $X$ into $\mathbb{P}_{n-1}$. Call the image of this embedding $C_n$. Let $S$ denote the open subspace of the Grassmannian variety consisting of $(n-3)$-dimensional linear subspaces of $\mathbb{P}_{n-1}$ not intersecting $C_n$. Namba showed that $S\times X$ is an $n^2$-sheeted unbranched covering space of $R_n/M$  \cite[Lemma 1.4.1]{Namba}. Note that for $n=3$, $S=\mathbb{P}_2\backslash C$. It seems possible to construct a manifold analogous to $N$. However, proving that this manifold is Oka is likely to be difficult.

Generalising our work to spaces of holomorphic maps $X\to \mathbb{P}$, where $X$ is a compact Riemann surface of genus $g\geq 2$ is another more ambitious project. Letting $R_n$ denote the set of degree $n$ maps $X\to\mathbb{P}$, Namba has shown that $R_n/M$ is biholomorphic to a specific open subset of a Grassmannian variety \cite[Corollary 1.3.13]{Namba}. As in the case of higher degree, constructing a manifold analogous to $N$ seems possible, but it is not at all clear whether it will be Oka, particularly because $X$ is no longer Oka. It is not known whether $X$ being hyperbolic is an obstruction to $R_n$ being Oka.


There are several smaller problems that arise from our work, but are not directly related to the study of holomorphic flexibility properties.

Firstly, it would be interesting to determine the covering transformations of the map $\Phi: W \to R_3$ to see if the branched covering map is normal. Our work shows that this boils down to showing that the covering $\lambda: N\to \mathbb{P}_2\backslash C$ is normal, so determining the covering transformations of this map is of interest.

Secondly, let $X_1$ and $X_2$ be two compact Riemann surfaces of genus $1$ and let $R_3(X_1)$, $R_3(X_2)$ be the associated complex manifolds of degree $3$ maps $X_1, X_2\to \mathbb{P}$. It seems unlikely that $R_3(X_1)$ and $R_3(X_2)$ are isomorphic unless $X_1$ is isomorphic to $X_2$. It would be nice to be able to confirm this.

Thirdly, we show that the set of maps in $R_3$ with double critical points forms a singular hypersurface of $R_3$, whose singular points are the maps with two double critical points. If $X$ is biholomorphic to $\mathbb{C}/\Gamma_0$, then the singular points of this hypersurface consists of maps that have three double critical points. To determine the isomorphism classes of these special maps is an interesting problem related to the theory of dessins d'enfants.

Finally, computing the fundamental group of $R_3$ seems possible because of the fibre bundle structure established by Namba. The key ingredient will be computing the fundamental group of $\mathbb{P}_2\backslash C$.
\chapter{Proof of Main Theorem\label{ch:two}}

\section{Topology}
Let $X,Y$ be topological spaces and $\mathcal{C}(X,Y)$ the set of all continuous maps $X\to Y$. The \emph{compact-open topology} on $\mathcal{C}(X,Y)$ is the smallest topology containing the sets
\[
\Set{f\in \mathcal{C}(X,Y)}{f(K)\subset U},
\]
where $K\subset X$ is compact and $U\subset Y$ is open.

The following results are well known \cite[p. 530]{Hatcher}.\\

\e{prop}\label{Metro}
If $X$ is compact and $Y$ is a metric space with metric $d$, then the compact-open topology on $\mathcal{C}(X,Y)$ is the metric topology given by the metric $d(f,g)=\sup\limits_{x\in X} d(f(x),g(x))$, the topology of uniform convergence.\\
\ee{prop}

\e{prop}
If $X$ is locally compact and $Z$ is any topological space, then a map $g:Z\to \mathcal{C}(X,Y)$ is continuous if and only if the map $f: Z\times X\to Y$, $f(z,x)=g(z)(x)$, is continuous.\\
\ee{prop}

Taking $g=\text{id}_{\mathcal{C}(X,Y)}$ in the above proposition implies that the evaluation map $e_\mathcal{C}:\mathcal{C}(X,Y)\times X\to Y$, $(f,x)\mapsto f(x)$, is continuous whenever $X$ is locally compact.

Let $X$ be a compact Riemann surface, $\mathbb{P}$ the Riemann sphere and $\mathcal{O}(X,\mathbb{P})$ the set of all holomorphic maps $X\to \mathbb{P}$. Then $\mathcal{O}(X,\mathbb{P})$ is a closed subspace of $\mathcal{C}(X,\mathbb{P})$, whose topology has a subbasis consisting of the sets
\[
[K,U]=\Set{f\in \mathcal{O}(X,\mathbb{P})}{f(K)\subset U},
\]
where $K\subset X$ is compact and $U\subset \mathbb{P}$ is open.

Let $R_n$ denote the subset of $\mathcal{O}(X,\mathbb{P})$ consisting of all maps of degree $n$. If the genus of $X$ is at least $1$, $R_1=\varnothing$, since a map $X\to \mathbb{P}$ of degree $1$ would be a biholomorphism. For $n\geq2$, $R_n$ may or may not be empty. In particular, $R_2\neq\varnothing$ if and only if $X$ is hyperelliptic. In general we can write
\[
\mathcal{O}(X,\mathbb{P})=R_0\sqcup R_2\sqcup R_3\sqcup \dots.
\]
Note that $R_0$ is just the set of constant maps. The following result is due to Namba \cite[Lemma 1.1.1]{Namba}. The proof, however, is our own work.\\

\e{prop}[Namba]\label{Open}
Let $X$ be a compact Riemann surface and $R_n$ denote the space of holomorphic maps $X\to\mathbb{P}$ of degree $n$. For each $n\geq0$, the set $R_n$ is both open and closed in $\mathcal{O}(X,\mathbb{P})$.
\ee{prop}
\e{proof}
It is sufficient to prove that $R_n$ is open. The set $R_0=[X,\mathbb{C}]\cup [X,\mathbb{P}\backslash \set{0}]$ is open. If $f\in R_n$, $n\geq 1$, then $f$ has $m$ distinct zeroes $p_1,\dots,p_m\in X$, where $1\leq m\leq n$. For $i=1,\dots,m$, let $U_i$ be a relatively compact neighbourhood (we take a neighbourhood to be open by definition) of $p_i$ in $X$ such that:
\begin{enumerate}
\item $\bar{U}_1,\dots , \bar{U}_m$ are pairwise disjoint,
\item each $U_i$ is a disc in some local coordinates,
\item each $\bar{U_i}$ contains no poles of $f$.
\end{enumerate}
Let $D=\bigcup\limits_{i=1}^m \bar{U_i}$. Then $D$ is compact and $\partial D= \bigcup\limits_{i=1}^m \partial{U_i}$. Let $\varepsilon=\min\limits_{\partial D}|f|>0$ and let
\[
W=\Set{g\in [D,\mathbb{C}]}{\max_{\partial D}|g-f|<\varepsilon}\cap [X\backslash \mathring{D}, \mathbb{P}\backslash\set{0}].
\]
It follows immediately that if $g\in W$, then $g$ has no zeroes on $\partial D$. Furthermore, the conditions 1--3 above allow us to apply Rouch\'e's theorem on each $\partial U_i$ to conclude that in each $\bar{U}_i$, every $g\in W$ has the same number of zeroes as $f$, counted with multiplicity. Since $g$ has no zeroes outside of $D$, it follows that $g\in R_n$. Hence $W\subset R_n$.

We shall now show that $W$ contains a neighbourhood of $f$. Since $f(D)$ is a compact subset of $\mathbb{C}$, there exist points $y_1,\dots,y_k\in D$ such that the open discs $B_{\frac{\varepsilon}{3}}(f(y_1))$,\dots, $B_{\frac{\varepsilon}{3}}(f(y_k))$ cover $f(D)$. Letting $K_i$ be the closure of $f^{-1}(B_{\frac{\varepsilon}{3}}(f(y_i)))$, we see that each $K_i$ is compact, $D\subset\Uns{i=1}{k}K_i$ and $f(K_i)\subset B_{\frac{\varepsilon}{2}}(f(y_i))$. Then setting $V_i=B_{\frac{\varepsilon}{2}}(f(y_i))$ and $S=\Inscts{i=1}{k}[K_i,V_i]\cap [X\backslash \mathring{D}, \mathbb{P}\backslash\set{0}]$ we see that $f\in S$. Since $S$ is open, it is sufficient to show that $S\subset W$. If $g\in S$ and $y\in K_i$, then since $g(y)\in V_i$,
\[
|g(y)-f(y_i)|<\frac{\varepsilon}{2}.
\]
Since $f(y)\in V_i$,
\[
|f(y)-f(y_i)|<\frac{\varepsilon}{2},
\]
so
\[
|f(y)-g(y)|<\varepsilon,
\]
implying that $g\in W$.
\ee{proof}

\section{Complex Structure}
\subsection{The Universal Complex Structure}
The following theorem is a special case of a deep result by Douady \cite[Section 10.2]{Douady}. A proof of this special case belongs to Kaup \cite[Theorem 1]{Kaup}. This theorem is fundamental to our work.\\

\e{thm}[Douady]\label{Douady}
Let $X$ be a compact complex manifold, let $Y$ be a reduced complex space and endow $\mathcal{O}(X,Y)$ with the compact-open topology. Then $\mathcal{O}(X,Y)$ can be given the structure of a reduced complex space satisfying the following universal property.
\begin{itemize}
\item[\emph{$(*)$}] If $Z$ is any reduced complex space, then a map $g:Z\to \mathcal{O}(X,Y)$ is holomorphic if and only if the map $f: Z\times X\to Y$, $f(z,x)=g(z)(x)$, is holomorphic.
\end{itemize}
\ee{thm}
Any complex structure on $\mathcal{O}(X,Y)$ satisfying $(*)$ is called a \emph{universal complex structure} and is unique up to biholomorphism.

Instead of using property $(*)$, Kaup defines a universal complex structure by the following two properties.
\begin{enumerate}
\item[(1)] The evaluation map $e:\mathcal{O}(X,Y)\times X\to Y$, $(f,x)\mapsto f(x)$, is holomorphic.
\item[(2)] If $Z$ is any reduced complex space, then a map $g:Z\to \mathcal{O}(X,Y)$ is holomorphic if the map $f: Z\times X\to Y$, $f(z,x)=g(z)(x)$, is holomorphic.
\end{enumerate}

Kaup's definition of a universal complex structure is equivalent to $(*)$. To prove this, firstly let $Z=\mathcal{O}(X,Y)$ and $g=\text{id}_Z$. Then $(*)$ implies that the evaluation map $e$ is holomorphic and so $(*)$ implies both (1) and (2). Conversely, in order for (1) and (2) to imply $(*)$ it must be shown that they imply the converse of (2). For $Z$ a reduced complex space and $g:Z\to \mathcal{O}(X,Y)$ holomorphic, the map $f: Z\times X\to Y$, $f(z,x)=g(z)(x)$, is the map $e\circ(g\times \text{id}_X)$, a composition of holomorphic maps.

Throughout this thesis we shall verify $(*)$ by verifying Kaup's properties, as the latter are almost always easier to check.

Property $(*)$ allows the following lemma to be proved.\\

\e{lem}\label{biholo}
Let $X,Y$ be complex manifolds with $X$ compact and let $\Aut(Y)$ denote the group of holomorphic automorphisms of $Y$. Then if $g\in \Aut(Y)$ is fixed, the map $\mathcal{O}(X,Y)\to \mathcal{O}(X,Y)$, $f\mapsto g\circ f$, is a biholomorphism.
\ee{lem}
\e{proof}
It is sufficient to show that the map is holomorphic. By property $(*)$ in Theorem $\ref{Douady}$, this holds if the map $\mathcal{O}(X,Y)\times X\to Y$, $(f,z)\mapsto (g\circ f)(z)$, is holomorphic. But this is just the map $g\circ e$, which is a composition of holomorphic maps.
\ee{proof}

Taking $X$ to be a Riemann surface of genus $g\geq 1$ and $Y=\mathbb{P}$, it follows immediately from Theorem \ref{Douady} and the previous section that
\[
\mathcal{O}(X,\mathbb{P})=R_0\sqcup R_2\sqcup R_3\sqcup \dots,
\]
possesses a universal complex structure. Furthermore, $R_n$ possesses a universal complex structure in the sense that it satisfies property $(*)$ with $R_n$ substituted for $\mathcal{O}(X,Y)$. Namba \cite[Chapter 1]{Namba} studied the universal complex structure on $R_n$ and proved the following theorem. See \cite[Proposition 1.1.4]{Namba}.\\

\e{thm}[Namba] \label{Manifold}
If $n\geq g$, then $R_n$ is a complex manifold of dimension $2n+1-g$.
\ee{thm}

For $n<g$, $R_n$ may be singular. Namba provides an explicit example in \cite[Section 2.2]{Namba}. The main goal of this section is to prove Theorem \ref{Manifold} when $X$ has genus $1$ (Theorem \ref{SCase}). Our proof of this special case is simpler than Namba's original proof, which covers the general case.

We shall finish this subsection by considering the complex structure of some basic examples of mapping spaces.

For our first example, we shall show that the natural identification of $R_0$ with $\mathbb{P}$ gives $R_0$ the universal complex structure. To show this, fix $x_0\in X$ and identify $R_0$ with $R_0\times\set{x_0}$. The map $R_0\times\set{x_0}\to\mathbb{P}$, $(f,x_0)\mapsto f(x_0)$, is obviously a bijection and is holomorphic when $R_0$ is given the universal complex structure. Kaup's property $(2)$ implies that the inverse map $\mathbb{P}\to R_0\times\set{x_0}$ is holomorphic if the map $\mathbb{P}\times X\to \mathbb{P}$, $(z, x)\mapsto z$, is holomorphic, which is obviously true.

In a similar vein, we shall show that the usual complex structure on the M\"obius group $M$ (the group of holomorphic automorphisms of $\mathbb{P}$) is the universal complex structure. The usual complex structure on $M$ is obtained by identifying it with $\mathbb{P}\GL_2(\mathbb{C})$. However, $M$ also is the set of degree $1$ holomorphic maps $\mathbb{P}\to\mathbb{P}$, so comes equipped with Douady's universal complex structure. To show that these two complex structures are the same, it suffices to check Kaup's properties $(1)$ and $(2)$ when $X=Y=\mathbb{P}$ and $\mathbb{P}\GL_2(\mathbb{C})$ is substituted for $\mathcal{O}(X,Y)$.
 
Firstly, we will verify property $(1)$. Let $e: \mathbb{P}\GL_2(\mathbb{C})\times\mathbb{P}\to \mathbb{P}$, $(f,z)\mapsto f(z)$, and $\tilde{e}:\GL_2(\mathbb{C})\times\mathbb{C}^2\backslash\set{0}\to\mathbb{C}^2\backslash\set{0}$, $(f,z)\mapsto f(z)$, be the evaluation maps and let $\Pi: \GL_2(\mathbb{C})\times \mathbb{C}^2\backslash\set{0}\to \mathbb{P}\GL_2(\mathbb{C})\times \mathbb{P}$ and $\pi: \mathbb{C}^2\backslash \set{0}\to \mathbb{P}$ be the projectivisation maps. Then the following diagram commutes.

\[
\begin{tikzcd}
\GL_2(\mathbb{C})\times \mathbb{C}^2\backslash\set{0} \arrow{r}{\tilde{e}} \arrow{d}{\Pi} & \mathbb{C}^2\backslash\set{0} \arrow{d}{\pi} \\
\mathbb{P}\GL_2(\mathbb{C})\times \mathbb{P} \arrow{r}{e} & \mathbb{P}
\end{tikzcd}
\]

It is obvious that the map $\tilde{e}$ is holomorphic. Then if $\sigma$ is a local holomorphic section of $\Pi$, it is clear that locally $e$ is the composition $\pi\circ \tilde{e}\circ \sigma$, which is holomorphic.

To show property $(2)$, let $Z$ be a reduced complex space and $g:Z\to \mathbb{P}\GL_2(\mathbb{C})$ such that $f: Z\times X\to\mathbb{P}$, $f(z,x)=g(z)(x)$, is holomorphic. Using a local holomorphic section of the projectivisation map $\GL_2(\mathbb{C})\to \mathbb{P}\GL_2(\mathbb{C})$, we can find a local lifting of $g$ to a map $\tilde{g}:U\to\GL_2(\mathbb{C})$, where $U$ is a member of a sufficiently fine open cover of $Z$. Then let $\tilde{f}:U\times \mathbb{C}^2\backslash\set{0}\to \mathbb{C}^2\backslash\set{0}$, $\tilde{f}(z,x)=\tilde{g}(z)(x)$, and let $\pi: \mathbb{C}^2\backslash\set{0}\to \mathbb{P}$ be the projectivisation map. Then the following diagram commutes.

\[
\begin{tikzcd}
U\times \mathbb{C}^2\backslash\set{0} \arrow{r}{\tilde{f}} \arrow{d}{\text{id}_U\times \pi} & \mathbb{C}^2\backslash\set{0} \arrow{d}{\pi} \\
U\times \mathbb{P} \arrow{r}{f} & \mathbb{P}
\end{tikzcd}
\]

Since $f$ is holomorphic, by taking local holomorphic sections of $\pi$ we see that $\tilde{f}$ is holomorphic. Now, the map $\tilde{g}$ is matrix-valued, so is holomorphic if and only if its component functions are holomorphic. Taking $x=(1,0)$, $(0,1)$, the holomorphicity of the map $\tilde{f}(\cdot,x): U\to\mathbb{C}^2\backslash\set{0}$ implies that these component functions are holomorphic. Then $\tilde{g}$ is holomorphic, so $g$ is holomorphic on $U$.

\subsection{Symmetric Products}\label{Sprod}
This subsection is drawn from \cite[Section 3.a]{Gunning} and \cite[Example 49 A.17 i]{KnK}.

Let $X$ be a Riemann surface and let $S_n$ denote the symmetric group of permutations of the set $\set{1,\dots, n}$. Then $S_n$ acts on the $n$-fold product $X^n$ by
\[
\sigma(z_1,\dots, z_n)=(z_{\sigma(z_1)}, \dots, z_{\sigma(z_n)}), \,\sigma\in S_n.
\]
The \emph{$n$-fold symmetric product} of $X$ is the orbit space $S^nX=X^n/S_n$. Then $S^n X$ is naturally identified with the set of positive divisors on $X$ of degree $n$. Since $S_n$ is a finite group, it follows that $S^nX$ is a complex space such that the quotient map $\pi:X^n\to S^nX$ is holomorphic \cite[Proposition 49A.16]{KnK}. The universal property of the quotient is that if $Y$ is any complex space, a holomorphic map $X^n\to Y$ factors as $\pi$ followed by a holomorphic map $S^n X\to Y$ if and only if it is $S_n$-invariant.

Note that $\pi$ is an open map. Indeed, if $V\subset X^n$ is open, then $\pi(V)$ is open if and only if the saturation of $V$ by $S_n$ is open. Since $S_n$ acts on $X^n$ via homeomorphisms, the saturation of $V$ is a union of open sets, hence open.


We shall now construct explicit charts on $S^n X$. They will be used later. They show that $S^nX$ is in fact a complex manifold. Consider the divisor $D=k_1p_1+\dots +k_mp_m\in S^nX$, where $p_1,\dots,p_m\in X$ are distinct and $k_1,\dots, k_m> 0$, $k_1+\dots+k_m=n$. Let $U_i$ be a neighbourhood of $p_i$ in $X$ such that $U_i\cap U_j = \varnothing$ whenever $i\neq j$. Furthermore, assume there are local coordinates $\varphi_i:U_i\to\mathbb{C}$. Let $U=U_1^{k_1}\times \dots\times U_m^{k_m}$. Then $\pi(U)$ is a neighbourhood of $D$ in $S^nX$.

If $\sigma\in S_n$, either $\sigma(U)=U$ or $\sigma(U)\cap U=\varnothing$. The set $G$ of all permutations $\sigma\in S_n$ such that $\sigma(U)=U$ is a subgroup of $S_n$. Then $\pi(U)=U/G$. 

Define $\rho_{i}: \mathbb{C}^{k_i}\to \mathbb{C}^{k_i}$, $\rho_{i}(z)=\left(\Sms{l=1}{k_i}z_l^1,\dots, \Sms{l=1}{k_i}z_l^{k_i}\right)$. Then let
\[
\varphi=\left( \rho_1\circ\underbrace{(\varphi_1\times \dots \times \varphi_1)}_{\text{$k_1$ times}}\right)\times\dots\times \left(\rho_m\circ\underbrace{(\varphi_m\times \dots \times \varphi_m)}_{\text{$k_m$ times}}\right).
\]
Then $\varphi$ is a holomorphic map $U\to\mathbb{C}^n$ and is clearly $G$-invariant, so induces a holomorphic map $\tilde{\varphi}:\pi(U)\to\mathbb{C}^n$. We wish to show that the image of $\tilde{\varphi}$ is open and $\tilde{\varphi}$ is a biholomorphism onto the image.

Firstly, we shall show that $S^n X$ is of pure dimension $n$. Let $\Omega\subset X^n$ be the subset of $(x_1,\dots, x_n)$ in $X^n$ such that for $i,j=1,\dots, n$, $x_i\neq x_j$ whenever $i\neq j$. Then $\Omega$ is open and dense in $X^n$. Furthermore, $S_n$ acts freely (and, since $S_n$ is finite, properly discontinuously) on $\Omega$, implying that $\Omega/S_n\subset S^n X$ is an $n$-dimensional connected complex manifold. Then the image of $\tilde{\varphi}$ is open and $\tilde{\varphi}$ is a biholomorphism onto the image if it is injective \cite[Proposition 46A.1]{KnK}.

We shall show that $\tilde{\varphi}$ is injective by showing for $i=1,\dots, m$ that if $\rho_i(z)=\rho_i(z')$, then $z'$ is a permutation of $z$. For $j=1,\dots, k_i$, let 
\[
s_{j}: \mathbb{C}^{k_i}\to \mathbb{C}, \quad s_j(z_1,\dots, z_{k_i})=\Sm{1\leq l_1<l_2<\dots<l_j\leq k_i}z_{l_1}z_{l_2}\cdots z_{l_j}.
\]
Then $s_j$ is the \emph{$j$th elementary symmetric function in $k_i$ variables}. It is well known that the polynomial
\begin{multline*}
(x-z_1)(x-z_2)\cdots (x-z_{k_i})\\=x^{k_i}-s_1(z_1,\dots, z_{k_i})x^{k_i-1}+\dots+ (-1)^{k_i}s_{k_i}(z_1,\dots, z_{k_i})
\end{multline*}
shows that the values $z_1,\dots, z_{k_i}$ are determined uniquely up to permutation by the values $s_1(z_1,\dots, z_{k_i})$, \dots, $s_{k_i}(z_1,\dots, z_{k_i})$. One of Newton's identities \cite{Newton} states that for $1\leq r\leq k_i$,
\[
\Sms{l=1}{k_i}z_l^r+\Sms{j=1}{r-1}\left((-1)^j s_j \Sms{l=1}{k_i}z_l^{r-j}\right)+(-1)^rrs_r=0.
\]
It follows that $s_1,\dots, s_{k_i}$ are uniquely determined by $\Sms{l=1}{k_i}z_l^1,\dots, \Sms{l=1}{k_i}z_l^{k_i}$, so $z_1,\dots, z_{k_i}$ are determined up to permutation by $\Sms{l=1}{k_i}z_l^1,\dots, \Sms{l=1}{k_i}z_l^{k_i}$. It follows that $\tilde{\varphi}$ is injective and hence is a chart on $S^nX$.

\subsection{The Jacobi Map}
Now assume that $X$ is a compact Riemann surface of genus $1$. Then there is a lattice $\Gamma$ in $\mathbb{C}$ such that $X$ is biholomorphic to $\mathbb{C}/\Gamma$. This allows us to identify $X$ with the complex Lie group $\mathbb{C}/\Gamma$. The Lie group structure on $X$ thus obtained is auxiliary and none of the results in this thesis depend on the choice of structure. For points $p_1,\dots,p_n\in X$ we can consider the divisor $p_1+\dots + p_n\in S^nX$, where here the symbol $+$ denotes the formal addition of points. Alternatively, we can consider the point $p_1+\dots + p_n\in X$, where $+$ denotes the Lie group addition in $X$. We trust that using the same symbol $+$ for both will not cause confusion. The \emph{Jacobi map} is the map
\[
\psi: S^nX\to X,\, k_1p_1+\dots + k_mp_m\mapsto k_1p_1+\dots+ k_m p_m,
\]
where $\psi$ maps the divisor $k_1p_1+\dots + k_mp_m$ to the Lie group sum $k_1p_1+\dots+ k_m p_m$. The map $\psi$ is holomorphic, as it is induced by the map $X^n\to X$, $(z_1,\dots,z_n)\mapsto z_1+\dots+z_n$, which is holomorphic and $S_n$-invariant.

It is well known that the map $\psi$ is a $\mathbb{P}_{n-1}$-bundle over $X$. We have been unable to track down a good reference for this fact, so will present our own argument.

We shall firstly show that $\psi$ is a bundle. Let $s_{n}:X\to X$, $z\mapsto nz$. Then $s_{n}$ is a holomorphic $n^2$-sheeted unbranched covering map. Fix $p\in X$ and let $U$ be an evenly covered neighbourhood of $p$. Then there exists a local section $\sigma: U\to X$ of $s_{n}$. Let
\[
h: U\times \psi^{-1}(0)\to\psi^{-1}(U), \,(x, p_1+\dots+p_n)\mapsto (p_1+\sigma(x))+\dots+(p_n+\sigma(x)).
\]
We claim that $h$ gives a local trivialisation of $\psi$.

To show that $h$ commutes with the projection $U\times \psi^{-1}(0)\to U$, note that for all $(x, p_1+\dots+p_n)\in U\times \psi^{-1}(0)$, we have $\psi(h(x, p_1+\dots+p_n))=n\sigma(x)=x$.

The map $h$ is clearly bijective, with inverse
\e{eqnarray*}
&h^{-1}: \psi^{-1}(U)\to U\times \psi^{-1}(0),\quad d_1+\dots+d_n\mapsto \biggl(\psi(d_1+\dots+d_n), \\&(d_1-\sigma(\psi(d_1+\dots+d_n)))+\dots+(d_n-\sigma(\psi(d_1+\dots+d_n)))\biggr).
\ee{eqnarray*}
Since the map $U\times X^n\to X^n$, $(x, p_1,\dots, p_n)\mapsto (p_1+\sigma(x), \dots, p_n+\sigma(x))$, is holomorphic, postcomposing it with the projection map $X^n\to S^nX$, we obtain an $S_n$-invariant holomorphic map $X^n\to S^nX$. This induces the map $h$, which is holomorphic by the universal property of the quotient. A similar argument shows that $h^{-1}$ is holomorphic and so is a trivialisation of $\psi$. Hence $\psi$ is a bundle over $X$. Finally, we determine its fibre.

Fix $p\in X$ and $D\in \psi^{-1}(p)$. Then since $\text{deg}(D)=n> 0$, we have $H^1(X, \mathcal{O}_D)=0$. By the Riemann-Roch theorem, $\text{dim}\, H^0(X, \mathcal{O}_D)=n$, so $\mathbb{P}H^0(X, \mathcal{O}_D)$ is biholomorphic to $\mathbb{P}_{n-1}$. Let $\set{\eta_1,\dots,\eta_n}$ be a basis of $H^0(X, \mathcal{O}_D)$. The following proposition is well known. See \cite[V.4 Proposition 4.3]{Griffy}.\\

\e{lem}\label{PJ}
The map 
\[
\kappa: \mathbb{P}_{n-1}\to \psi^{-1}(p),\, [a_1, \dots, a_n]\mapsto (a_1\eta_1+\dots+a_n\eta_n) +D,
\]
is a biholomorphism.
\ee{lem}

It follows immediately that $\psi$ is a $\mathbb{P}_{n-1}$-bundle over $X$. Since both $\mathbb{P}_{n-1}$ and $X$ are Oka, $\psi$ is a fibre bundle with Oka fibre and Oka base. So we have the following.\\

\e{prop}
Let $X$ be a compact Riemann surface of genus $1$. Then for $n\geq1$, $S^n X$ is an Oka manifold.
\ee{prop}

For $f\in R_n$, let $D_0(f)=\text{max}\set{0, (f)}\in S^n X$ be the divisor of zeroes of $f$, and let $D_\infty(f)=\text{max}\set{0, -(f)}$ be the divisor of poles of $f$.

The following lemma is part of the proof of Theorem 1.2.1 in \cite{Namba}. We present a simplified version of Namba's proof.\\

\e{lem}[Namba]\label{Divisors}
The maps $D_0, D_\infty: R_n\to S^nX$ are holomorphic.
\end{lem}
\e{proof}
We shall show that $D_0$ is holomorphic. Fix $f\in R_n$ and let $D_0(f)=k_1p_1+\dots+ k_m p_m$, where $k_1+\dots+k_m=n$ and $p_1,\dots, p_m\in X$ are distinct. Then, as in the proof of Proposition \ref{Open}, we can take a neighbourhood $W$ of $f$ in $R_n$, such that there exist pairwise disjoint neighbourhoods $U_1,\dots, U_m$ of $p_1,\dots, p_m$ respectively, such that every $g\in W$ has the same number of zeroes as $f$ in each $U_i$, and $g$ has no poles in any of the sets $\bar{U}_1,\dots, \bar{U}_m$. Furthermore, we can assume that in local coordinates, $U_1,\dots, U_m$ are discs.

Since $R_n$ is a reduced complex space, we can take $W$ so small that it can be identified with a subvariety of an open subset $\Omega$ of $\mathbb{C}^k$ for some $k$. Shrinking $\Omega$ if necessary, we can ensure that for every $i=1,\dots, m$, the evaluation map $e_i: W\times U_i \to \mathbb{C}$ extends to a holomorphic map
\[
\tilde{e}_i: \Omega\times U_i\to\mathbb{C}.
\]
For $y\in \Omega$, let $\tilde{e}_{i,y}: U_i\to \mathbb{C}$, $\tilde{e}_{i,y}(w)=\tilde{e}_i(y,w)$. It follows from the proof of Proposition \ref{Open} that we may assume that $\Omega$ is so small that if $\varepsilon=\min\limits_{\partial U_i}|f|>0$, then for all $y\in \Omega$,
\[
\max\limits_{\partial U_i}|\tilde{e}_{i,y}-f|<\epsilon.
\]
Then Rouch\'e's theorem implies that $\tilde{e}_{i,y}$ has $k_i$ zeroes in $U_i$. Let these zeroes in local coordinates be denoted $z_{i1}(y),\dots, z_{ik_i}(y)$, where $\Sm{i}k_i=m$. It is a well-known application of the argument principle \cite[p.~153--154]{Al} that
\[
\Sms{j=1}{k_i}z_{ij}(y)^p=\frac{1}{2\pi i}\igrl{\partial U_i}\rous{\tilde{e}_i}{w}(y,w)\cdot \frac{w^p}{\tilde{e}_i (y,w)}\,dw
\]
for $p=1,\dots, k_i$. So $\Sms{j=1}{k_i}z_{ij}(y)^p$ defines a holomorphic function $\Omega\to\mathbb{C}$. But if $y\in W$, then
\[
\left(\Sms{j=1}{k_1}z_{1j}(y),\dots, \Sms{j=1}{k_1}z_{1j}(y)^{k_1}, \dots, \Sms{j=1}{k_m}z_{mj}(y),\dots, \Sms{j=1}{k_m}z_{mj}(y)^{k_m}\right)
\]
are the local coordinates of $D_0(y)$ under the standard chart on $S^n X$ centred at $D_0(f)$ (see Subsection \ref{Sprod}). Hence $D_0$ is holomorphic. The map $D_\infty$ can be shown to be holomorphic by applying the above argument to the reciprocal of the evaluation map.
\ee{proof}

\subsection{The Theta Function}

A reference for this subsection is \cite[p.~34--35, p.~50]{Mir}.

Let $\Gamma$ be a lattice in $\mathbb{C}$. We may assume that $\Gamma=\mathbb{Z} + \tau \mathbb{Z}$, for $\tau\in \mathbb{C}$ with $\text{Im}(\tau)>0$. Then let
\[
\theta(z)=\Sms{n=-\infty}{\infty}e^{\pi i [n^2\tau + 2n z]}.
\]
This series converges absolutely and uniformly on compact subsets of $\mathbb{C}$ and so defines
an entire function $\theta: \mathbb{C}\to\mathbb{C}$, the \emph{Jacobi theta function of $\Gamma$}.

The function $\theta$ is classical with well-understood properties. For every $z\in\mathbb{C}$,
\[
\theta(z+1)=\theta(z)\, \quad \text{and} \quad\,\theta(z+\tau)= e^{-\pi i (\tau +2z)}\theta(z).
\]
Furthermore, $\theta$ has a simple zero at $\frac{1}{2}(1+\tau)$. The translation properties of $\theta$ imply that $\theta$ has a simple zero at each point on the translated lattice $\frac{1}{2}(1+\tau)+\Gamma$. These are the only zeroes of $\theta$.

For $x\in\mathbb{C}$, let $\theta^{(x)}(z)=\theta(z-\frac{1}{2}(1+\tau)-x)$. Then $\theta^{(x)}$ has simple zeroes on $x+\Gamma$.

The following theorem relates theta functions with elliptic functions. For the proof, see \cite[p.~34--35, p.~50]{Mir}.

\e{thm}\label{Class}
If $x_1, \dots, x_n, y_1,\dots, y_n\in\mathbb{C}$ such that $\Sms{i=1}{n}(x_i-y_i)\in \Gamma$, then
\[
f(z)=\Pds{i=1}{n}\theta^{(x_i)}(z)\bigg/\Pds{i=1}{n}\theta^{(y_i)}(z),
\]
is a $\Gamma$-periodic meromorphic function on $\mathbb{C}$. Conversely, suppose that $g$ is a $\Gamma$-periodic meromorphic function on $\mathbb{C}$ with zeroes $x_1,\dots,x_n$ modulo $\Gamma$ and poles $y_1,\dots,y_n$ modulo $\Gamma$, counted with multiplicity. Then $\Sms{i=1}{n}(x_i-y_i)\in \Gamma$ and there exists $c\in \mathbb{C}^*$ such that $g=cf$, where $f$ is the meromorphic function defined above.
\ee{thm}

The above theorem is essentially an explicit version of Abel's theorem for the torus \cite[Section 20.8]{Forster}. Indeed, let $X$ be a compact Riemann surface of genus $1$. Recall the zero and polar divisor maps $D_0, D_\infty:R_n\to S^nX$ and the Jacobi map $\psi: S^nX\to X$. Then Abel's theorem for the torus states that if $D_1,D_2\in S^nX$ have no common points, then there exists a holomorphic map $f:X\to \mathbb{P}$ with $D_0(f)=D_1$ and $D_\infty(f)=D_2$ if and only if $\psi(D_1)=\psi(D_2)$.

Let 
\[
B_n=\Set{(D_1,D_2)\in S^nX\times S^n X}{\text{$D_1$ and $D_2$ have no points in common}}
\]
for $n\geq2$. Then $B_n$ is an open subset of $S^nX\times S^n X$. Let
\[
C_n=\Set{(D_1,D_2)\in S^nX\times S^n X}{\psi(D_1)=\psi(D_2)}.
\]
Then $C_n$ is a $(2n-1)$-dimensional subvariety of $S^nX\times S^nX$. Also, $Q_n= B_n\cap C_n$ is a $(2n-1)$-dimensional subvariety of $B_n$.

By Lemma \ref{Divisors} the map $(D_0, D_\infty): R_n\to Q_n$ is holomorphic. By Abel's theorem, its image lies in $Q_n$. We will now provide a holomorphic section of this map on a neighbourhood of each point of $Q_n$. This is needed for Theorem \ref{Divisor} below.

Let $D_1=p_1+\dots+ p_n$ and $D_2= q_{1} +\dots+ q_{n}$ be divisors on $X$ such that $(D_1,D_2)\in Q_n$. Then we can take neighbourhoods $U_i$, $V_i$ around each $p_i$, $q_i$ respectively, such that $U_i\cap V_j=\varnothing$, $U_i=U_j$ if $p_i=p_j$, $V_i=V_j$ if $q_i=q_j$, $U_i\cap U_j=\varnothing$ if $p_i\neq p_j$ and $V_i\cap V_j=\varnothing$ if $q_i\neq q_j$. Furthermore, we may assume that there exist holomorphic sections $\xi_i:U_i\to \mathbb{C}$ and $\eta_i:V_i\to \mathbb{C}$ of the quotient map $\pi: \mathbb{C}\to \mathbb{C}/\Gamma$.

If $(p_1',\dots, p_n')\in U_1\times\dots\times U_n$ and $(q_1',\dots, q_n')\in V_1\times\dots\times V_n$, such that $(p_1'+\dots+p_n', q_1'+\dots+q_n')\in Q_n$, then $\Sms{i=1}{n}(\xi_i(p_i)-\eta_i(q_i))\in \Gamma$. We can define a function $\Theta: (U_1\times\dots\times U_n\times V_1\times\dots\times V_n)\cap \pi^{-1}(Q_n)\to R_n$ by
\[
\Theta(p_1',\dots, p_n', q_1',\dots, q_n')= \Pds{i=1}{n}\theta^{(\xi_i(p_i'))}\bigg/\Pds{i=1}{n}\theta^{(\eta_i(q_i'))},
\]
where $\pi: X^n\times X^n\to S^nX\times S^nX$ is the projection. The map $\Theta$ is invariant under the action of $S_n\times S_n$, so it induces a local section of the map $(D_0, D_\infty):R_n\to Q_n$.

By the universal property $(*)$ in Theorem \ref{Douady}, the local section is holomorphic if the map $\Theta$ is holomorphic. But the map $\Theta$ is holomorphic if the map 
\begin{align*}
&(U_1\times\dots\times U_n\times V_1\times\dots\times V_n)\times X\to\mathbb{P},\\&((p_1,\dots, p_n, q_1,\dots, q_n), z)\mapsto \Theta(p_1,\dots, p_n, q_1,\dots, q_n)(z),
\end{align*}
is holomorphic. But this follows from the fact that
\[
\theta^{(\xi_i(p_i'))}(z)=\theta(z-\frac{1}{2}(1+\tau)-\xi_i(p_i')),
\]
which varies holomorphically with $p_i'$ and $z$. Similarly, $\theta^{(\eta_i(q_i'))}(z)$ varies holomorphically with $q_i'$ and $z$.

\subsection{The Divisor Map}

The following theorem and its proof can be found in the work of Namba \cite[proof of Theorem 1.2.1]{Namba}. Our proof is simpler than Namba's proof, as it only covers the special case when $X$ has genus $1$.\\

\e{thm}\label{Divisor}
The map $(D_0, D_\infty):R_n\to Q_n$, $n\geq2$, is a holomorphic principal $\mathbb{C}^*$-bundle.
\ee{thm}
\e{proof}
Lemma $\ref{Divisors}$ implies that the map $(D_0,D_\infty)$ is holomorphic. Given any pair of divisors $(D_1,D_2)\in Q_n$, there exists a neighbourhood $U$ of $(D_1,D_2)$ in $Q_n$ with a holomorphic section $\Theta: U\to R_n$ defined as above.

Let $h: \mathbb{C}^*\times U\to R_n$,
\[
h(\alpha, (D_1',D_2'))=\alpha\Theta(D_1',D_2').
\]
Theorem \ref{Class} implies that $h$ is a bijection $\mathbb{C}^*\times U\to (D_0,D_\infty)^{-1}(U)$. We claim that $h$ is holomorphic. Let $\Phi_1: \mathbb{C}^* \times R_n \to R_n$, $(\alpha, g)\mapsto \alpha g$. Then $h=\Phi_1\circ (\text{id}_{\mathbb{C}^*}\times \Theta)$, so $h$ is holomorphic if $\Phi_1$ is. By the universal property $(*)$ in Theorem \ref{Douady}, $\Phi_1$ is holomorphic if the map $\Phi_2:\mathbb{C}^*\times R_n\times X\to \mathbb{P}$, $(\alpha, g, x)\mapsto \alpha g(x),$ is holomorphic. But $\Phi_2=\Phi_3\circ (\text{id}_{\mathbb{C}^*}\times e)$, where $e$ is the evaluation map and $\Phi_3: \mathbb{C}^*\times \mathbb{P}\mapsto \mathbb{P}$, $(\alpha, z)\mapsto \alpha z$. Since $\Phi_3$ is clearly holomorphic, it follows that $h$ is as well.

We now show that $h^{-1}$ is holomorphic and hence is a local trivialisation of $(D_0,D_\infty)$. For any $f\in (D_0,D_\infty)^{-1}(U)$, ${\Theta(D_0(f), D_\infty (f))}$ is a meromorphic function with the same zeroes and poles as $f$. Hence $f/\Theta(D_0(f), D_\infty (f))$ is a nonzero constant function. We can then define a map
\[
\kappa: (D_0,D_\infty)^{-1}(U)\to \mathbb{C}^*,\, \kappa(f)= f /\Theta(D_0(f), D_\infty (f)).
\]
So $h^{-1}=(\kappa, (D_0,D_\infty))$.

It remains to show that $\kappa$ is holomorphic. Fixing $f\in (D_0, D_\infty)^{-1}(U)$, choose $a\in X$ such that $f(a)\in\mathbb{C}^*$. Take the neighbourhood
\[
W=\Set{g\in R_n}{g(a)\in\mathbb{C}^*}\cap (D_0, D_\infty)^{-1}(U)
\]
of $f$. On $W$, $\kappa(g)=g(a)/{\Theta(D_0(g), D_\infty (g))(a)}$ is holomorphic in $g$.
\ee{proof}

\subsection{Nonsingularity of $R_n$}

As before, let $X$ be a compact Riemann surface of genus $1$ and for $n\geq 2$, let $R_n$ denote the set of degree $n$ holomorphic maps $X\to\mathbb{P}$. We can now prove the following theorem, a special case of Theorem \ref{Manifold}. See \cite[Proposition 1.1.4]{Namba}.\\

\e{thm}[Namba]\label{SCase}
For $n\geq 2$, the complex space $R_n$ is a complex manifold of dimension $2n$.
\ee{thm}
\e{proof}
Choose $f\in R_n$. We shall show that $f$ has a smooth neighbourhood. Since $f$ has finitely many branch points, Lemma \ref{biholo} allows us to assume that $D_0(f)$ and $D_\infty(f)$ each consist of $n$ distinct points. Theorem \ref{Divisor} implies that it is sufficient to prove that $Q_n$ is nonsingular in a neighbourhood of the point $(D_0(f), D_\infty(f))$. Let $D_0(f)=p_1+\dots+p_n$ and $D_\infty(f)=q_1+\dots+q_n$ where $p_i\neq p_j$ and $q_i\neq q_j$ when $i\neq j$. Take pairwise disjoint neighbourhoods $U_1,\dots, U_n$ and $V_1,\dots, V_n$ of $p_1,\dots, p_n$ and $q_1,\dots, q_n$ respectively such that there are local coordinates, $x_i: U_i \to \mathbb{C}$ and $y_i: V_i \to \mathbb{C}$ at $p_i$ and $q_i$ respectively. Then we know from Subsection \ref{Sprod} that the map $x_1\times \dots \times x_n\times y_1\times \dots \times y_n$ induces a chart on $S^nX\times S^nX$ at $(D_1, D_2)$. This chart provides a biholomorphism of $Q_n$ with the set of points $(x_1,\dots, x_n, y_1,\dots, y_n)$ in $U_1\times \dots\times U_n \times V_1\times \dots\times V_n$ such that $\Sms{i=1}{n}x_i-y_i\in \Gamma$. This set is the intersection of a union of disjoint hyperplanes in $\mathbb{C}^{2n}$ with an open subset and so is smooth. It follows that $Q_n$ is a manifold of dimension $2n-1$. Then by Theorem \ref{Divisor}, $R_n$ is nonsingular of dimension $2n$.
\ee{proof}

\subsection{The Space $R_n/M$}

This subsection is drawn from the work of Namba \cite[Section 1.3]{Namba}. Let $R_n$ be as before, with $n\geq2$, and let $M$ be the group of M\"obius transformations. Recall that $M$ acts on $R_n$ by postcomposition. Clearly the action is free. Our goal in this subsection is to show that the orbit space $R_n/M$ is a complex manifold such that the projection map $\pi: R_n\to R_n/ M$ is a principal $M$-bundle. For this we shall use the following special case of a theorem of Holmann \cite[Satz 21, 24]{Holmann}.\\

\e{thm}[Holmann]\label{Holmann}
Let $X$ be a complex space and let $G$ be a complex Lie group acting freely and properly on $X$. Then the orbit space $X/G$ possesses a complex structure such that the quotient map $\pi: X\to X/G$ is a holomorphic principal $G$-bundle. Furthermore, if $X$ is a complex manifold, then so is $X/G$.
\ee{thm}

Note that since a principal $G$-bundle always possesses local holomorphic sections, it follows that $X/G$ has the universal property of the quotient: if $Y$ is any complex space, a holomorphic map $X\to Y$ factors as $\pi$ followed by a holomorphic map $X /G\to Y$ if and only if it is $G$-invariant.

Let $X$ and $Y$ be compact complex manifolds. Let $O(X,Y)$ denote the set of all open holomorphic maps $X\to Y$. Then $O(X,Y)$ is an open subset of $\mathcal{O}(X,Y)$ \cite[Proposition 1]{NOpen} and so is naturally equipped with Douady's complex structure as in Theorem \ref{Douady}. Furthermore, it is well known that the group $\Aut(Y)$ is a complex Lie group and acts holomorphically on $O(X,Y)$ by postcomposition \cite[Theorem D]{Chu}. Namba \cite[Theorem 1]{NOpen} used the previous theorem to prove the following result.\\

\e{thm}[Namba]\label{AutY}
The orbit space $O(X,Y)/\Aut(Y)$ admits a complex structure such that the quotient map $\pi: O(X,Y)\to O(X,Y)/\Aut(Y)$ is a holomorphic principal $\Aut(Y)$-bundle.
\ee{thm} 

We shall prove Theorem \ref{AutY} in our special case, when $X$ is a compact Riemann surface and $Y=\mathbb{P}$. This allows our proof to be more elementary than Namba's proof. Then $O(X,Y)$ is just the set of nonconstant holomorphic maps $X\to Y$. Furthermore, $\Aut(Y)$ is the group $M$ of M\"obius transformations. Lemma \ref{biholo} implies that $M$ acts holomorphically on $O(X,Y)$. By Theorem \ref{Holmann} it remains to show the map
\[
\Phi: M\times O(X,\mathbb{P})\to O(X,\mathbb{P})\times O(X,\mathbb{P}), (g,f)\mapsto (g\circ f, f),
\]
is proper. Let $K\subset O(X,\mathbb{P})\times O(X,\mathbb{P})$ be compact. Let $(g_k,f_k)$ be a sequence in $\Phi^{-1}(K)$. Since $K$ is compact, we may assume that $(g_k\circ f_k, f_k)$ converges to $(h, f)$ in $K$. Furthermore, Proposition \ref{Open} allows us to assume that the $f_k$, $f$ and $h$ are all of the same degree $n\geq 1$.





We can identify $M$ with $\mathbb{P}\text{GL}_2(\mathbb{C})$. If $M_n(\mathbb{C})$ denotes the vector space of $n\times n$ matrices with complex entries, there is a natural inclusion map $\iota: \mathbb{P}\text{GL}_2(\mathbb{C})\to \mathbb{P}\text{M}_2(\mathbb{C})$. Since $\mathbb{P}\text{M}_2(\mathbb{C})$ is compact, by passing to a subsequence, we can assume there exists $A\in \mathbb{P}\text{M}_2(\mathbb{C})$ such that $g_{k}\to A$ in $\mathbb{P}\text{M}_2(\mathbb{C})$. If $\text{det}(A)\neq 0$, then $A\in M$ and so $g_{k}\to A$ in $M$ and we are done.

Otherwise $\text{det}(A)=0$. Then $\text{ker}(A)$ corresponds to a unique point $z_0\in \mathbb{P}$. We claim that for some $L\in \mathbb{P}$, $g_{k}(z)\to L$ for every $z\in \mathbb{P}\backslash\set{z_0}$. Indeed, we claim that $g_{k}\to L$ uniformly on compact subsets of $\mathbb{P}\backslash\set{z_0}$. To show this, we can assume that $z_0=\infty$. Take $a_{k}, b_{k}, c_{k}, d_{k}\in\mathbb{C}$ with $a_{k}d_{k}-b_{k}c_{k}\neq0$ such that for all $z\in\mathbb{P}$, $g_{k}(z)=(a_{k}z+b_{k})/(c_{k}z+d_{k})$. We can assume that $|a_k|^2+|b_k|^2+|c_k|^2+|d_k|^2=1$. Then, passing to a subsequence, we can assume that $a_{k}, b_k, c_k,d_k$ converge to $a,b,c,d\in \mathbb{C}$ respectively, where $|a|^2+|b|^2+|c|^2+|d|^2=1$. Since $z_0=\infty$, it follows that $a=c=0$ and so $|b|^2+|d|^2=1$. Hence $L=b/d$.

If $d\neq0$, let $\tilde{K}\subset\mathbb{C}$ be compact and let $C=\sup\limits_{z\in \tilde{K}}|z|$. Then for all $z\in \tilde{K}$,
\e{eqnarray*}
\left|{\frac{a_kz+b_k}{c_kz+d_k}-\frac{b}{d}}\right|&=&\left|\frac{(da_k-bc_k)z+(b_kd-bd_k)}{d(c_kz+d_k)}\right|\\
&\leq& \frac{|da_k-bc_k|\cdot C}{|d|\cdot|c_kz+d_k|}+\frac{|b_kd-bd_k|}{|d|\cdot|c_kz+d_k|}.
\ee{eqnarray*}
Let $\varepsilon>0$. Then for $k$ large enough, $|da_k-bc_k|<\varepsilon$ and $|b_kd-bd_k|<\varepsilon$. Furthermore, since $|c_kz|\leq |c_k|\cdot C$, for $k$ large enough, $|c_kz|<|d|/4\neq0$ and $|d_k|>3|d|/4$. Then
\e{eqnarray*}
|c_kz+d_k|\geq \left||c_kz|-|d_k|\right|>\frac{|d|}{2},
\ee{eqnarray*}
and so $|c_kz+d_k|^{-1}< 2/{|d|}$. Then
\[
\left|{\frac{a_kz+b_k}{c_kz+d_k}-\frac{b}{d}}\right|\leq \frac{2C\varepsilon}{|d|^2}+\frac{2\varepsilon}{|d|^2}.
\]

This shows that (after passing to a subsequence) the sequence $(g_k)$ converges to the constant map $L=b/a$ uniformly on compact subsets of $\mathbb{C}=\mathbb{P}\backslash \set{z_0}$.

If $d=0$, then the above argument shows that $1/g_k$ converges uniformly on compact subsets of $\mathbb{C}$ to $0$.
 
It follows that $h$ is constant, contradicting $h$ having degree $n\geq 1$.

Since the action of $M$ on $O(X,\mathbb{P})$ maps $R_n$ to $R_n$, it follows that $R_n/M$ is a complex manifold such that the projection map $\pi: R_n\to R_n/M$ is a holomorphic principal $M$-bundle. Furthermore, $R_n/M$ possesses the universal property of the quotient.

\section{Elliptic Functions of Degree 2}

\subsection{Elliptic Curves}

Let $X$ be a compact Riemann surface of genus $1$. As before, let $R_n$ denote the set of degree $n$ holomorphic maps $X\to\mathbb{P}$. The goal of this section is to show that $R_2$ is an Oka manifold. In fact, we will show that $R_2$ is a homogeneous space.

For this section we will represent $X$ as the quotient of $\mathbb{C}$ by a lattice $\Gamma$. We will begin by reviewing some basics.

Let $\omega_1,\omega_2\in\mathbb{C}$ be linearly independent over $\mathbb{R}$. Then $\Gamma=\omega_1\mathbb{Z}+\omega_2\mathbb{Z}$ is a lattice in $\mathbb{C}$. We say that $\Gamma$ is \emph{square} if $\omega_1=e^{i\pi/2}\omega_2$ and that $\Gamma$ is \emph{hexagonal} if $\omega_1=e^{i\pi/3}\omega_2$. 

The quotient $\mathbb{C}/\Gamma$ is a Riemann surface of genus $1$, a $1$-dimensional complex torus. Furthermore, $\mathbb{C}/\Gamma$ is a complex Lie group and the quotient map $\pi: \mathbb{C}\to \mathbb{C}/\Gamma$ is holomorphic, a Lie group homomorphism and the universal covering map.

For two lattices $\Gamma, \Gamma'$ in $\mathbb{C}$, every non-constant holomorphic map $f:\mathbb{C}/\Gamma\to\mathbb{C}/\Gamma'$ fixing $0$ lifts to a holomorphic map of the form $\tilde{f}:\mathbb{C}\to\mathbb{C}$, $z\mapsto \alpha z$, where $\alpha\in\mathbb{C}^*$ and $\alpha\Gamma\subset \Gamma'$. It follows that $f$ is a homomorphism of Lie groups. Furthermore, $f$ is a biholomorphism if and only if $\alpha\Gamma=\Gamma'$. Let $A_\Gamma$ denote the group of automorphisms of $\mathbb{C}/\Gamma$ that fix $0$. The following result is proved in \cite[p.~64]{Mir}.\\

\e{prop}\label{hex}
If $f\in A_\Gamma$ lifts to $z\mapsto \alpha z$, then:
\e{enumerate}
\item $\alpha\in\set{\pm1,\pm i}$ if $\Gamma$ is a square lattice,
\item $\alpha\in\set{\pm1,e^{\pm i\frac{\pi}{3}},e^{\pm i\frac{2\pi}{3}}}$ if $\Gamma$ is a hexagonal lattice,
\item $\alpha\in\set{\pm1}$ otherwise.
\ee{enumerate}
Thus $A_\Gamma$ contains $4$ elements if $\Gamma$ is a square lattice, $6$ elements if $\Gamma$ is hexagonal and $2$ elements otherwise.
\ee{prop}

If $t\in \mathbb{C}/\Gamma$, let $\tau_t:\mathbb{C}/\Gamma\to\mathbb{C}/\Gamma$, $z\mapsto z+t$, be the translation map. Then $\tau_t$ is a holomorphic automorphism of $\mathbb{C}/\Gamma$. Let $T_\Gamma =\Set{\tau_t}{t\in\mathbb{C}/\Gamma}$, a subgroup of the full automorphism group $\Aut(\mathbb{C}/\Gamma)$. Then $\Aut(\mathbb{C}/\Gamma)=T_\Gamma \rtimes A_\Gamma$. Note that $T_\Gamma$ is naturally identified with $\mathbb{C}/\Gamma$ itself; however we find it convenient to distinguish the two.

For the Riemann sphere $\mathbb{P}$, consider $\mathcal{O}(\mathbb{C}/\Gamma,\mathbb{P})$, the set of holomorphic functions $f:\mathbb{C}/\Gamma\to\mathbb{P}$. If $f\in \mathcal{O}(\mathbb{C}/\Gamma,\mathbb{P})$ is not equal to the constant function $\infty$, we say $f$ is an \emph{elliptic function with respect to} $\Gamma$. A meromorphic function $g:\mathbb{C}\to\mathbb{P}$ is said to be \emph{doubly periodic with respect to} $\Gamma$ if for all $\omega\in \Gamma$ and $z\in\mathbb{C}$, $f(z)=f(z+\omega)$. It is clear that the set of elliptic functions and doubly periodic meromorphic functions with respect to $\Gamma$ are in bijective correspondence via the equality $g=f\circ\pi$. The set of elliptic functions with respect to $\Gamma$ then forms a field under the usual multiplication and addition of meromorphic functions.

Recall that the automorphism group of $\mathbb{P}$ is $M$, the M\"obius group of transformations $g:\mathbb{P}\to\mathbb{P}$ given by the formula
\[
g(z)=\frac{az+b}{cz+d},
\]
for $a,b,c,d\in \mathbb{C}$ with $ad-bc\neq0$. A M\"obius transformation is determined uniquely by where it takes three distinct points. The group $M$ acts on $\mathcal{O}(\mathbb{C}/\Gamma,\mathbb{P})$ via postcomposition, with $g\cdot f=g\circ f$, $f\in \mathcal{O}(\mathbb{C}/\Gamma,\mathbb{P})$, $g\in M$. Similarly, $\text{Aut}(\mathbb{C}/\Gamma)$ acts on $\mathcal{O}(\mathbb{C}/\Gamma,\mathbb{P})$ via precomposition with $f\cdot \tau=f\circ \tau^{-1}$, for $\tau\in\text{Aut}(\mathbb{C}/\Gamma)$. The product group $M\times \text{Aut}(\mathbb{C}/\Gamma)$ then acts on $\mathcal{O}(\mathbb{C}/\Gamma,\mathbb{P})$ in the obvious way: $(g,\tau)\cdot f=g\circ f\circ\tau ^{-1}$. Since the elements of $M$ and $\text{Aut}(\mathbb{C}/\Gamma)$ are of degree $1$, it follows that $R_n$ is invariant under the action of $M\times \text{Aut}(\mathbb{C}/\Gamma)$.

\subsection{The Weierstrass $\wp$-function}
We shall use \cite{FB} as a general reference for this subsection.

Letting $\Gamma'=\Gamma\backslash\set{0}$ and $k\in\mathbb{Z}$, define the \emph{Eisenstein series} $G_k(\Gamma)$ by
\[
G_k(\Gamma)=\Sm{\omega\in \Gamma'}\omega^{-k}.
\]
The series converges absolutely for $k\geq3$.

The \emph{Weierstrass} $\wp$-\emph{function} is the map $\wp:\mathbb{C}\backslash\Gamma\to \mathbb{C}$ defined by
\[
\wp(z)=\frac{1}{z^2}+\Sm{\omega\in \Gamma'} \frac{1}{(z-\omega)^2}-\frac{1}{\omega^2}.
\]
This series converges uniformly on compact subsets of $\mathbb{C}\backslash\Gamma$ and thus is a meromorphic function on $\mathbb{C}$ with poles of order $2$ precisely on the lattice $\Gamma$. Furthermore, $\wp$ is an even function with $\wp(-z)=-\wp(z)$ for all $z\in \mathbb{C}$. Differentiating term by term, we see that
\[
\wp'(z)=-2\Sm{\omega\in \Gamma}\frac{1}{(z-\omega)^3}.
\]
It is clear that $\wp'$ is doubly periodic with respect to $\Gamma$ and an odd function, with $\wp'(-z)=-\wp'(z)$ for all $z\in \mathbb{C}$. Since $\wp$ is even, it follows that $\wp$ is also doubly periodic with respect to $\Gamma$.

The functions $\wp$, $\wp'$ induce elliptic functions and hence members of $\mathcal{O}(\mathbb{C}/\Gamma,\mathbb{P})$. The $\wp$-function is of degree $2$, while $\wp'$ is of degree $3$. These two functions are related by the following theorem.\\

\e{thm}
Let $g_2=60G_4$ and $g_3=140G_6$. Then $\wp$ and $\wp'$ satisfy the \emph{Weierstrass equation}
\[
(\wp')^2=4\wp^3-g_2\wp -g_3.\\
\] 
\ee{thm}

The functions $\wp$ and $\wp'$ generate the field of elliptic functions over $\mathbb{C}$.\\

\e{thm}\label{FieldThm}
The field of elliptic functions on $\mathbb{C}/\Gamma$ is $\mathbb{C}(\wp,\wp')$.\\
\ee{thm}

\subsection{Degree 2 Case}

Let $\omega_3=\omega_1+\omega_2\in \Gamma$, and for $i=1,2,3$, let $b_i=\pi(\frac{\omega_i}{2})\in\mathbb{C}/\Gamma$. Then if $z\in\mathbb{C}/\Gamma$ and $z=-z$, it follows that either $z=0$ or $z=b_i$ for some $i$. Now $\wp'$ is an odd function of degree $3$, implying that $\wp'$ has a simple zero at $b_1$, $b_2$, $b_3$. Then $0$, $b_1$, $b_2$, $b_3$ are branch points of $\wp$.

Recall the Hurwitz formula \cite[Section 17.14]{Forster} for an $n$-sheeted branched covering of the Riemann sphere by a genus $g$ surface:
\[
g=\frac{b}{2}-n+1,
\]
where $b$ is the branching order of the covering. In the case $g=1$ we obtain
\[
b=2n.
\]
Since $\wp$ is of degree $2$, this implies that $0$, $b_1$, $b_2$, $b_3$ are the only branch points of $\wp$.

Recall that the action of $M\times \text{Aut}(\mathbb{C}/\Gamma)$ on $\mathcal{O}(\mathbb{C}/\Gamma,\mathbb{P})$ leaves $R_d$ invariant. Since $\mathbb{C}/\Gamma$ is a subgroup of $\text{Aut}(\mathbb{C}/\Gamma)$, there is an induced action of $M\times \mathbb{C}/\Gamma$ on $R_2$. It is shown in \cite{Helicoid} that this action is transitive.\\

\e{prop}
The action of $M\times \mathbb{C}/\Gamma$ on $R_2$ is transitive.
\ee{prop}
We provide a short proof for the convenience of the reader.
\e{proof}
If $f\in R_2$, then there exists $(g_1,t)\in M\times \mathbb{C}/\Gamma$ such that $f_1=g_1\circ f\circ \tau_t^{-1}\in R_2$ has a double pole at $0$ and is zero at $b_1$. There also exists $g_2\in M$ such that $f_2=g_2\circ \wp\in R_2$ has a double pole at $0$ and a double zero at $b_1$. Then the quotient $f_1/f_2$ is an elliptic function with a pole of order at most $1$ at $b_1$ and is holomorphic elsewhere. Since $R_1=\varnothing$, it follows that $f_1=cf_2$ for some $c\in\mathbb{C}^*$. Taking $g=g_1^{-1}\circ c g_2$, we see that $f=g\circ \wp\circ \tau_{-t}^{-1}$ and so the action is transitive.
\ee{proof} 

We shall now determine the stabiliser of $\wp$ under the action of $M\times \mathbb{C}/\Gamma$. Recall that $\wp$ has a double pole at $0$. Let $y\in \mathbb{C}/\Gamma$ be a zero of $\wp$. Since $\wp$ is even, the other zero of $\wp$ is at $-y$ (which is equal to $y$ only if $\Gamma$ is the square lattice).

Now, if $(g,t)\in M\times \mathbb{C}/\Gamma$ such that
\begin{equation}\label{Stabby}
g\circ \wp \circ \tau_t^{-1}=\wp,
\end{equation}
then the set of branch points of $g\circ \wp \circ \tau_t^{-1}$ is $S=\set{0,b_1,b_2,b_3}$. This implies that $\tau_t$ acts as a permutation on $S$ and so $t\in S$. Hence $t=-t$, which implies that $\wp\circ\tau_t^{-1}$ is even, since
\[
\wp\circ\tau_t^{-1}(-z)=\wp(-z-t)=\wp(z+t)=\wp(z-t)=\wp\circ\tau_t^{-1}(z).
\]
Let $h\in M$ such that $h\circ \wp\circ \tau_t^{-1}(0)=\infty$ and $h\circ \wp\circ \tau_t^{-1}(y)=0$. Then $h\circ \wp\circ \tau_t^{-1}(-y)=0$ and so $f=h\circ \wp\circ \tau_t^{-1}$ has the same zeroes and poles as $\wp$. Hence $\wp/f$ is a constant $c\in \mathbb{C}^*$. Choosing $h$ to be the unique M\"obius transformation for which this constant is $1$, we have shown that for any $t\in S$, there exists a unique $g=g_t\in M$ satisfying (\ref{Stabby}). We have proved the following proposition.\\

\e{prop}
The stabiliser of $\wp$ under the action of $M\times \mathbb{C}/\Gamma$ on $R_2$ is
\[
\Stab(\wp)=\Set{(g_t,t)}{t\in S}.
\]
This group is isomorphic to $\mathbb{Z}_2\times \mathbb{Z}_2$.
\ee{prop}

Now equip the set $R_2$ with the universal complex structure as in Theorem \ref{Douady}. We wish to show that $R_2$ is biholomorphic to $(M\times \mathbb{C}/\Gamma)/ \text{Stab}(\wp)$.

Firstly, we shall show $M\times \mathbb{C}/\Gamma$ acts holomorphically on $R_2$. By property $(*)$ in Theorem \ref{Douady}, the map $M\times \mathbb{C}/\Gamma\times R_2\to R_2$, $(g, t, f)\mapsto g\circ f\circ\tau_{-t}$, is holomorphic if the map $M\times \mathbb{C}/\Gamma\times R_2\times \mathbb{C}/\Gamma\to \mathbb{P}$, $(g, t, f,x)\mapsto g(f(x-t))$, is holomorphic. But this follows from the fact that the evaluation map $e: R_2\times \mathbb{C}/\Gamma\to \mathbb{P}$ is holomorphic.

Now, consider the map $\Phi: M\times \mathbb{C}/\Gamma\to R_2$, $(g,t)\mapsto g\circ\wp\circ\tau_{-t}$. It is holomorphic and surjective. Furthermore, $\Phi$ is invariant under the action of $\Stab(\wp)$ and so descends to a holomorphic bijection $\tilde{\Phi}:(M\times \mathbb{C}/\Gamma)/\Stab(\wp)\to R_2$. Since $M\times \mathbb{C}/\Gamma$ is connected, $R_2$ is connected and so has pure dimension. It follows that $\tilde{\Phi}$ is a biholomorphism \cite[Proposition 49.13]{KnK}. Hence we have proved the following theorem.\\

\e{thm}
Let $X$ be a compact Riemann surface of genus $1$ and $R_2$ the complex manifold of degree $2$ maps $X\to\mathbb{P}$. Then $R_2$ is a homogeneous space and hence is an Oka manifold. Moreover, $R_2$ is biholomorphic to $(M\times X)/(\mathbb{Z}_2\times \mathbb{Z}_2)$.
\ee{thm}

See \cite[Proposition 5.5.1]{Forstneric} for the proof that every homogeneous space is Oka.

\section{A $9$-sheeted Covering of $R_3/ M$}\label{Rum}

The goal of this section is to show that there is a smooth cubic curve $C$ in $\mathbb{P}_2$ such that $\mathbb{P}_2\backslash C\times X$ is a $9$-sheeted unbranched covering space of $R_3/M$. This result is due to Namba and we draw upon his work \cite[Sections 1.3, 1.4]{Namba}.

In this chapter, it is more natural for us to think of $C$ as a curve in the dual projective plane $\mathbb{P}_2^\vee$, which is of course biholomorphic to $\mathbb{P}_2$. We shall let $[a,b,c]$ denote the homogeneous coordinates of the point in $\mathbb{P}_2^\vee$ corresponding to the line $ax_1+bx_2+cx_3=0$ in $\mathbb{P}_2$.

Let $X$ be a compact Riemann surface of genus $1$. Let $\psi: S^n X\to X$, $n\geq2$, be the Jacobi map and $D_0, D_\infty$ the maps $R_n\to S^n X$ taking a function to its divisor of zeroes and poles respectively. Let $\beta= \psi\circ D_0: R_n\to X$. Then $\beta$ is holomorphic and $M$-invariant, so by the universal property of the quotient, it induces a holomorphic map $\chi: R_n/M\to X$ such that the following diagram commutes.

\[
\begin{tikzcd}
R_n \arrow{r}{D_0} \arrow{dr}{\beta} \arrow[swap]{d}{\pi} & S^nX \arrow{d}{\psi} \\
R_n/M \arrow{r}{\chi} & X
\end{tikzcd}
\]

From now on, we shall take $n=3$. Fixing $p\in X$ and $D\in \psi^{-1}(p)$, since $X$ has genus $g=1$ and $\text{deg}(D)=3>2g-2=0$, it follows that $H^1(X, \mathcal{O}_D)=0$ \cite[Theorem 17.16]{Forster}. The Riemann-Roch theorem implies that $\text{dim}\, H^0(X, \mathcal{O}_D)=3$. Thus $\mathbb{P}H^0(X, \mathcal{O}_D)$ is biholomorphic to $\mathbb{P}_2^\vee$. Let $\set{\eta^1,\eta^2,\eta^3}$ be a basis of $H^0(X, \mathcal{O}_D)$. Then the map $\eta_p: X\to \mathbb{P}_2^\vee$, $x\mapsto [\eta^1(x), \eta^2(x), \eta^3(x)]$, is a holomorphic embedding of $X$ in $\mathbb{P}_2^\vee$ \cite[Theorem 17.22]{Forster}. We shall denote the image of $X$ under this map by $C_p$. By the degree-genus formula, $C_p$ is a cubic curve \cite[Section V.3, Proposition 2.15]{Mir}.\\

Now Lemma \ref{PJ} states that the map
\[
\kappa: \mathbb{P}_{2}\to \psi^{-1}(p),\quad [a_1, a_2, a_3]\mapsto (a_1\eta^1+a_2\eta^2+a_3\eta^3)+D,
\]
is a biholomorphism. Let $\nu=\kappa^{-1}$.\\

\e{prop}\label{Group Law}
Let $x_1, x_2, x_3\in \psi^{-1}(p)$. Then $x_1+x_2+x_3=p$ if and only if there is a line passing through $\eta_p(x_1), \eta_p(x_2), \eta_p(x_3)$ that intersects $C_p$ at $\eta_p(x_i)$ with multiplicity equal to the multiplicity of $x_i$ in $x_1+x_2+x_3$, for $i=1,2,3$.
\ee{prop}
\e{proof}
Let $x_1+x_2+x_3=p$. Then Abel's theorem implies that there exists $f\in \beta^{-1}(p)$ such that $D_0(f)=x_1+x_2+x_3$. We claim that the point $[a_1, a_2, a_3]=\nu(D_0(f))\in \mathbb{P}_2$ is dual to the desired line passing through $\eta_p(x_1), \eta_p(x_2), \eta_p(x_3)$. Indeed,
\[
D_0(f)= (a_1\eta^1+a_2\eta^2+a_3\eta^3)+D.
\]
Now if $D$ does not contain the point $x_i$, then $a_1\eta^1(x_i)+a_2\eta^2(x_i)+a_3\eta^3(x_i)=0$ with multiplicity equal to the multiplicity of $x_i$ in $D_0(f)$ and we are done. Otherwise, we can assume that $D$ contains the point $x_1$ with multiplicity $n$, and $D_0(f)$ contains $x_1$ with multiplicity $m$, for $1\leq n,m\leq 3$. Then $a_1\eta^1+a_2\eta^2+a_3\eta^3$ has order $m-n$ at $x_1$. Now, we know that $H^0(X,\mathcal{O}_D)$ is globally generated as $X$ has genus $g=1$ and $\text{deg}(D)=3\geq2g=2$ \cite[Theorem 17.19]{Forster}. Hence one of $\eta^1, \eta^2, \eta^3$ will have a pole at $x_1$ of order $n$. Then in local coordinates $z$ on $X$ about $x_1$ we can find holomorphic functions $\tilde{\eta}^j$, for $j=1,2,3,$ such that locally $\eta^j(w)=\tilde{\eta}^j(w)/(w-z(x_1))^n$. Then $\eta_p(x_1)=[\tilde{\eta}^1(x_1), \tilde{\eta}^2(x_1), \tilde{\eta}^3(x_1)]$. We have
\[
a_1\tilde{\eta}^1(w)+a_2 \tilde{\eta}^2(w)+a_3\tilde{\eta}^3(w)= (w-z)^n(a_1\eta^1(w)+a_2\eta^2(w)+a_3\eta^3(w)).
\]
It follows that the above expression has order $m$ at $x_1$ and so we are done.

Conversely, suppose $[a_1, a_2, a_3]$ defines a line in $\mathbb{P}_2^\vee$ which intersects $C_p$ at the points $\eta_p(x_1)$, $\eta_p(x_2)$, $\eta_p(x_3)$ with the desired multiplicities, where $x_1, x_2, x_3\in X$ are taken appropriately. Then for $i=1,2,3$, it is easy to use the computations above to show that $(a_1\eta^1+a_2\eta^2+a_3\eta^3)+D=x_1+x_2+x_3$. Then $x_1+x_2+x_3\in \psi^{-1}(p)$ and we are done.
\ee{proof}

Now, any cubic curve in $\mathbb{P}_2^\vee$ can be given a group law by choosing an inflection point to be the identity $0$ and then defining all collinear triples of points to sum to $0$. Proposition \ref{Group Law} shows that when $p=0$, the group structure induced by $X$ on $C_0$ by the embedding $\eta_0$ is such a group law. Indeed, $\eta_0$ necessarily maps $0\in X$ to an inflection point of $C_0$ as the divisor $0+0+0\in\psi^{-1}(0)$ determines a line intersecting $C_0$ at $\eta_0(0)$ with multiplicity $3$.

We shall now construct a $9$-sheeted unbranched covering map $\mathbb{P}_2^\vee\backslash C_0\times X\to R_3/M$. The construction relies on the following result of Namba \cite[Lemma 1.3.6]{Namba}.\\

\e{prop}[Namba]\label{Cubic Complement}
For all $p\in X$, there exists a biholomorphism $\chi^{-1}(p)\to \mathbb{P}_2^\vee\backslash C_p$.
\ee{prop}

Namba's original proof covered the case when $X$ was a Riemann surface of higher genus. We shall present our own proof of this result as it is significantly simpler when $X$ is of genus $1$. 
\e{proof}
Consider the following commuting diagram:

\[
\begin{tikzcd}
\beta^{-1}(p) \arrow {drr}{\mu} \arrow{r}{(D_0,D_\infty)} \arrow{d}{\pi} & (\psi^{-1}(p)\times \psi^{-1}(p))\backslash S \arrow{r}{\nu\times \nu} & (\mathbb{P}_2\times \mathbb{P}_2)\backslash \Delta \arrow{d}{L} \\
\chi^{-1}(p) \arrow{rr}{\tilde{\mu}} && \mathbb{P}_2^\vee
\end{tikzcd}
\]

If $x,y\in\mathbb{P}_2$ are distinct points, there exists a unique line $L(x,y)\in \mathbb{P}_2^{\vee}$ passing through $x$ and $y$. Furthermore, $L$ depends holomorphically on $x$ and $y$. Letting $\Delta$ denote the diagonal of $\mathbb{P}_2\times \mathbb{P}_2$, we obtain a holomorphic map $L: (\mathbb{P}_2\times \mathbb{P}_2)\backslash \Delta \to \mathbb{P}_2^\vee$.

Let $S$ denote the set of pairs of divisors in $\psi^{-1}(p)\times \psi^{-1}(p)$ that share at least one point. Then Abel's theorem implies that $(D_0, D_\infty)$ maps $\beta^{-1}(p)$ onto $(\psi^{-1}(p)\times \psi^{-1}(p))\backslash S$. Lemma \ref{PJ} implies that the map
\[
\kappa: \mathbb{P}_{2}\to \psi^{-1}(p),\quad [a_1, a_2, a_3]\mapsto (a_1\eta_1+a_2\eta_2+a_3\eta_3)+D,
\]
is a biholomorphism. Then taking $\nu=\kappa^{-1}$, we see that $\nu\times\nu$ maps $(\psi^{-1}(p)\times \psi^{-1}(p))\backslash S$ biholomorphically onto an open subset of $(\mathbb{P}_2\times \mathbb{P}_2)\backslash \Delta$ (it is not easy, nor necessary for us, to explicitly describe this open subset).

Let $\mu=L\circ (\nu\times \nu)\circ (D_0, D_\infty)|\beta^{-1}(p)$. Then $\mu$ is a holomorphic map $\beta^{-1}(p)\to \mathbb{P}_2^\vee$. We shall show that $\mu$ is $M$-invariant. Take $f\in \beta^{-1}(p)$ and $g\in M$. Let $r_0=D_0(g\circ f)$, $r_\infty=D_\infty(g\circ f)$ and let $L_f=\mu(f)$. Then it is sufficient to show that the points $\nu(r_0), \nu(r_\infty)$ are both on $L_f$.

We shall firstly show that $\nu(r_0)$ is on $L_f$. If $\nu(D_0(f))=[a_1, a_2, a_3]$ and $\nu(D_\infty(f))=[b_1, b_2, b_3]$, it is sufficient to show that for some $\lambda, \lambda'\in \mathbb{C}$, $\nu(r_0)=[\lambda a_1+\lambda'b_1, \lambda a_2+\lambda'b_2, \lambda a_3+\lambda'b_3]$. We can assume that $D=D_\infty(f)$, so $f\in H^0(X, \mathcal{O}_D)$. Then take $a,b,c,d\in \mathbb{C}$ with $ad-bc\neq0$ such that $g\circ f = (af+b)/(cf+d)$. Then if $a=0$, $r_0=D_\infty(f)$, and we are done. Otherwise $a\neq 0$ and for $z\in X$, $(g\circ f)(z)=0$ if $f(z)=-b/a$. Now, for some $h\in H^0(X, \mathcal{O}_D)$ we have $r_0=(h)+D_\infty(f)$. Since $r_0$ is disjoint from $D_\infty(f)$, then $D_\infty(f)=D_\infty(h)$ and $r_0=D_0(h)$. The function $af+b$ also has the same poles and zeroes as $h$, implying that for some $C\in \mathbb{C}^*$, $h=C(af+b)$. But
\[
af+b=a(a_1\eta_1+a_2\eta_2+a_3\eta_3)+b(b_1\eta_1+b_2\eta_2+b_3\eta_3),
\]
which implies that $\nu(r_0)=\nu(D_0(af+b))=[aa_1+bb_1, aa_2+bb_2, aa_3+bb_3]$, which is what we wanted to show. A similar argument using $D_\infty$ shows that $r_\infty$ is also on $L_f$.

Note that for each $p\in X$, $M$ acts on $\beta^{-1}(p)$, so $\beta^{-1}(p)/M=\chi^{-1}(p)$. Since $\mu$ is $M$-invariant, it induces a holomorphic map $\tilde{\mu}:\chi^{-1}(p)\to \mathbb{P}_2^\vee$.

We shall now prove that $\tilde{\mu}$ yields the desired biholomorphism in Proposition \ref{Cubic Complement}. Firstly, we shall show that $\tilde{\mu}(\chi^{-1}(p))=\mathbb{P}^\vee_2\backslash C_p$. This is equivalent to $\mu(\beta^{-1}(p))=\mathbb{P}_2^\vee\backslash C_p$. Take $(a,b)\in \mathbb{P}_2\times \mathbb{P}_2$ such that $L(a, b)\in C_p$ and let $a=[a_1, a_2, a_3]$, $b=[b_1, b_2, b_3]$. Then for some $z\in X$, $L(a,b)=[\eta_1(z), \eta_2(z), \eta_3(z)]$. Then $a_1\eta_1(z)+a_2\eta_2(z)+a_3\eta_3(z)=0$ and $b_1\eta_1(z)+b_2\eta_2(z)+b_3\eta_3(z)=0$. For fixed $D\in\psi^{-1}(p)$, $\kappa(a)=(a_1\eta_1+a_2\eta_2+a_3\eta_3)+D$ and $\kappa(b)=(b_1\eta_1+b_2\eta_2+b_3\eta_3)+D$. But the functions $a_1\eta_1+a_2\eta_2+a_3\eta_3$ and $b_1\eta_1+b_2\eta_2+b_3\eta_3$ in $H^0(X,\mathcal{O}_D)$ share the zero $z\in X$. Hence $\kappa(a)$ and $\kappa(b)$ share the point $z$ and so $(\kappa(a), \kappa(b))\in S$. Hence $\mu(\beta^{-1}(p))\subset\mathbb{P}_2^\vee\backslash C_p$.

Now to show that each $L(a,b)\in \mathbb{P}_2^\vee\backslash C_p$ is in $\mu(\beta^{-1}(p))$, Abel's theorem implies that it is sufficient to show that $\kappa(a)$ and $\kappa(b)$ are disjoint. Let $a= [a_1, a_2, a_3]$, $b=[b_1, b_2, b_3]$. Then it is equivalent to show that the functions $a_1\eta_1+a_2\eta_2+a_3\eta_3$ and $b_1\eta_1+b_2\eta_2+b_3\eta_3$ have no common zeroes. But if for some $z\in X$, $(a_1\eta_1+a_2\eta_2+a_3\eta_3)(z)=0=(b_1\eta_1+b_2\eta_2+b_3\eta_3)(z)$, then $L(a, b)=[\eta_1(z), \eta_2(z), \eta_3(z)]\in C_p$. Hence $\mu(\beta^{-1}(p))=\mathbb{P}_2^\vee\backslash C_p$.

To show that $\tilde{\mu}$ is injective, it is sufficient to show that if $x\in \mathbb{P}_2^\vee\backslash C_p$ and $f_1, f_2\in \mu^{-1}(x)$, then there exists $g\in M$ such that $f_1=g\circ f_2$. Let $\nu(D_0(f_1))=a=[a_1, a_2, a_3]$ and $\nu(D_\infty(f_1))=b=[b_1, b_2, b_3]$. We can assume that $D=D_\infty(f_1)$. Then $f_1\in H^0(X, \mathcal{O}_D)$ and $D_0(f_1)=(f_1)+D$. Then for some $C\in \mathbb{C}^*$, $a_1\eta_1+a_2\eta_2+a_3\eta_3= C f_1$ and $b_1\eta_1+b_2\eta_2+b_3\eta_3$ is a constant function. Then if $\nu(D_0(f_2))=[\alpha_1, \alpha_2, \alpha_3]$ and $\nu(D_\infty(f_2))=[\beta_1, \beta_2, \beta_3]$, there exist coefficients $c_1, c_2, c_3, c_4\in\mathbb{C}$ such that
\begin{align*}
\alpha_1\eta_1+\alpha_2\eta_2+\alpha_3\eta_3 &= c_1(a_1\eta_1+a_2\eta_2+a_3\eta_3)+c_2(b_1\eta_1+b_2\eta_2+b_3\eta_3)\\
&=c_1 f_1+ c_2
\end{align*}
and $\beta_1\eta_1+\beta_2\eta_2+\beta_3\eta_3 = c_3f_1+ c_4$. Since $D_0(f_2), D_\infty(f_2)$ are distinct, $c_1c_4-c_2c_3\neq0$, and $g(z)=(c_1z+c_2)/(c_3z+c_4)$ is a M\"obius transformation. Moreover it is easily seen that $g\circ f_1$ has the same zeroes and poles as $f_2$ and so there exists a M\"obius transformation taking $f_1$ to $f_2$.

To show that $\tilde{\mu}^{-1}$ is holomorphic, we shall show that for any $x\in\mathbb{P}_2^\vee$ there is a local holomorphic section of $\mu$ on some neighbourhood of $z$. It is clear that there is a neighbourhood $U$ of $x$ on which a local holomorphic section $\sigma_1$ of $L$ can be defined. As in the proof of Theorem $\ref{Divisor}$, there also exists a local holomorphic section $\sigma_2$ of the map $(D_0, D_\infty)$ in a neighbourhood of any point in $(\psi^{-1}(p)\times \psi^{-1}(p))\backslash S$. Shrinking $U$ if necessary, we see that $\sigma_2\circ (\kappa\times \kappa)\times \sigma_1$ is the desired local holomorphic section of $\mu$.
\ee{proof}

We shall continue to take $\mu: \beta^{-1}(p)\to \mathbb{P}_2^\vee$ to be as in the proof of Proposition \ref{Cubic Complement}. The following proposition reveals the relationship between the geometry of $\mathbb{P}_2^\vee\backslash C_p$ and $\chi^{-1}(p)$.\\

\e{prop}\label{LineLemma}
Let $q\in \mathbb{P}_2^\vee\backslash C_p$ and $x_1, x_2, x_3\in X$ such that $x_1+x_2+x_3=p$. Then $q$ is on the line that intersects $C_p$ at the points $\eta_p(x_1), \eta_p(x_2), \eta_p(x_3)$ counted with multiplicity if and only if $\set{x_1, x_2, x_3}$ is a fibre of every map in $\mu^{-1}(q)$.
\ee{prop}
\e{proof}
Suppose that $[a_1, a_2, a_3]$ defines a line through $\eta_p(x_1)$, $\eta_p(x_2)$, $\eta_p(x_3)$ and $q$ with the right multiplicity. Then we know from the proof of Proposition \ref{Group Law} that $\kappa[a_1, a_2, a_3]=x_1+x_2+x_3$. Now $q\in\mathbb{P}_2^\vee\backslash C_p$ is a line in $\mathbb{P}_2$ passing through $[a_1, a_2, a_3]$. If $[b_1, b_2, b_3]$ is a point on $q$, distinct from $[a_1, a_2, a_3]$, then we claim that $\kappa[b_1, b_2, b_3]$ is a divisor distinct from $x_1+x_2+x_3$. Indeed, if $\kappa[a_1, a_2, a_3]$ and $\kappa [b_1, b_2, b_3]$ shared $x_1$, say, then the corresponding lines in $\mathbb{P}_2^\vee$ would intersect at $\eta_p(x_1)\in C_p$. Since these lines already intersect at $q$, it follows that they are the same line, which contradicts the fact that they are distinct. Abel's theorem implies that there is a meromorphic function $f\in \beta^{-1}(p)$ such that $D_0(f)=\kappa [a_1, a_2, a_3]=x_1+x_2+x_3$ and $D_\infty(f)=\kappa [b_1, b_2, b_3]$. Indeed it is clear from the definition of $\mu$ that $\mu(f)=q$. Moreover, in the proof of Proposition \ref{Cubic Complement} we showed that $\mu^{-1}(q)$ is the $M$-orbit of $f$. Hence $\set{x_1, x_2, x_3}$ is a fibre of every $f\in \mu^{-1}(q)$.

Conversely suppose that $q$ is such that $\set{x_1, x_2, x_3}$ is a fibre of every $f\in \mu^{-1}(q)$. Then we can choose $f\in \mu^{-1}(q)$ such that $D_0(f)=x_1+x_2+x_3$. Let $[a_1, a_2, a_3]= \nu(x_1+x_2+x_3)$. Then the definition of $\mu$ implies that $q$ is a line in $\mathbb{P}_2$ passing through $[a_1, a_2, a_3]$. The dual statement is that the line defined by $[a_1, a_2, a_3]$ passes through $q$ in $\mathbb{P}_2^\vee$. The proof of Proposition \ref{Group Law} implies that this line intersects $C_p$ at the points $\eta_p(x_1)$, $\eta_p(x_2)$, $\eta_p(x_3)$ with the right multiplicity and we are done.
\ee{proof}

Having proved Proposition \ref{Cubic Complement}, we will now construct a covering map $\mathbb{P}_2^\vee\backslash C_0 \times X\to R_3/M$. Let $s_3: X\to X, x\mapsto 3x$. Then $s_3$ is a holomorphic $9$-sheeted unbranched covering map. The following lemma is due to Namba \cite[Lemma 1.4.1]{Namba}. Our proof is simpler and more straightforward than Namba's original proof.\\

\e{prop}[Namba]\label{Top Cover}
There exists a map $\tilde{s}_3: \mathbb{P}_2^\vee\backslash C_0\times X\to R_3/ M$, which makes the following diagram commute:
\[
\begin{tikzcd}
\mathbb{P}_2^\vee\backslash C_0 \times X \arrow{r}{\tilde{s}_3} \arrow{d}{\proj_2} & R_3/M \arrow{d}{\chi} \\
X \arrow{r}{s_3} & X
\end{tikzcd}
\]
where $\proj_2$ is the projection onto the second factor. Furthermore, this diagram is a pullback square. 
\ee{prop}
\e{proof}
By Proposition \ref{Cubic Complement} we can identify $\mathbb{P}_2^\vee\backslash C_0$ with $\chi^{-1}(0)$. For notational convenience we will suppress this isomorphism. Let $P$ denote the pullback of $R_3/M$ along $s_3$. Then it is enough to show that $P$ is biholomorphic to $\chi^{-1}(0)\times X$. By definition,
\[
P=\Set{(\psi,x)\in (R_3/M)\times X}{\chi(\psi)=s_3(x)}.
\]
If $\pi: R_3\to R_3/M$ is the projection map, then since $\pi$ is an $M$-bundle, about any $\psi\in R_3/M$ we can find a neighbourhood $U$ on which there exists a local holomorphic section $\sigma: U\to R_3$ of $\pi$. For $p\in X$, let $\tau_p:X\to X$, $x\mapsto x+p$. Then we can define the map
\[
\Phi: \chi^{-1}(0)\times X\to P, \, (\psi, x)\mapsto (\pi(\sigma(\psi)\circ\tau_{-x}),x).
\]
To show that $\Phi$ is well defined, we must show that it is independent of the choice of the section $\sigma$ and show that the image of $\chi^{-1}(0)\times X$ is contained in $P$. Take $\psi\in \chi^{-1}(0)$ and two local holomorphic sections $\sigma_1, \sigma_2$ of $\pi$ defined on some neighbourhood $U$ of $\psi$. Then if $\varphi\in U$, there exists $g\in M$ such that $\sigma_2(\varphi)=g\circ\sigma_1(\varphi)$. Then
\[
\pi(\sigma_2(\varphi)\circ\tau_{-x})=\pi(g\circ\sigma_1(\varphi)\circ \tau_x)=\pi(\sigma_1(\varphi)\circ \tau_{-x}),
\]
so the map is independent of the choice of section. To show that the image is contained in $P$, note that if the zeroes of $\sigma(\psi)$ are $p_1$, $p_2$, $p_3$, then the zeroes of $\sigma(\psi)\circ\tau_{-x}$ are $p_1+x$, $p_2+x$, $p_3+x$. Then $\chi(\sigma(\psi)\circ\tau_{-x})=3x=s_3(x)$ and so $\Phi$ is well defined.

The map $\Phi$ is clearly a bijection with inverse
\[
\Phi^{-1}: P\to \chi^{-1}(0)\times X,\, (\psi,x)\mapsto (\pi(\sigma(\psi)\circ\tau_{x}),x).
\]
Both $\Phi$ and $\Phi^{-1}$ are obviously holomorphic and so $P$ is biholomorphic to $\chi^{-1}(0)\times X$.
\ee{proof}

If we identify $\mathbb{P}_2^\vee\backslash C_0$ with $\chi^{-1}(0)$ as in the proof above, then every element of $\mathbb{P}_2^\vee\backslash C_0$ can be written as $\pi(f)$, where $f\in R_3$ is chosen appropriately. Unravelling the construction of $\tilde{s}_3$ in the proof above, we get $\tilde{s}_3(\pi(f),x)=\pi(f\circ\tau_{-x})$. This formula will be used later.

Since $\tilde{s}_3$ is surjective and $\mathbb{P}_2^\vee\backslash {C_0}$, $X$ and $M$ are connected, we can finally conclude that $R_3$ is connected.\\

\e{cor}
The complex manifold $R_3$ is connected.
\ee{cor}

Proposition \ref{Top Cover} also implies that $R_3$ is Oka if and only if $\mathbb{P}_2^\vee\backslash C_0$ is Oka. Indeed $M$ is Oka and $\pi: R_3\to R_3/M$ is a principal $M$-bundle. Then $R_3$ is Oka if and only if $R_3/M$ is Oka. Since $\tilde{s}_3$ is an unbranched covering map, $R_3/M$ is Oka if and only $\mathbb{P}_2^\vee\backslash C_0\times X$ is Oka. Since $X$ is Oka, $\mathbb{P}_2^\vee\backslash C_0\times X$ is Oka if and only if $\mathbb{P}_2^\vee\backslash C_0$ is Oka.\\

\e{cor}
The complex manifold $R_3$ is Oka if and only if the cubic complement $\mathbb{P}_2^\vee\backslash C_0$ is Oka.
\ee{cor}

Whether the complement of a smooth cubic curve in $\mathbb{P}_2$ is Oka or not is a well-known open problem in Oka theory \cite[Open Problem B, p. 20]{Survey}. Hence, we are unable to conclude whether $R_3$ is Oka or not.

From now on, we will blur the distinction between $\mathbb{P}_2$ and $\mathbb{P}_2^\vee$.

\section{Oka Branched Covering Space of $\mathbb{P}_2\backslash C$}\label{OBC}

Let $F:\mathbb{C}^3\to \mathbb{C}$ be a homogeneous polynomial of degree $3$. Then the equation 
\e{equation}\label{Cubic}
F(x,y,z)=0
\ee{equation}
defines a projective curve $C$ in $\mathbb{P}_2$. We shall assume that the partial derivatives $F_x,F_y, F_z$ have no common zeroes on $C$, so $C$ is smooth. Then $C$ is biholomorphic to a complex torus $\mathbb{C}/\Gamma$, for some lattice $\Gamma$ in $\mathbb{C}$. We may assume that $F$ is written in the standard form
\begin{equation}\label{Standard}
F(x,y,z)=y^2z-4x^3+g_2xz^2+g_3z^3,
\end{equation}
where $g_2, g_3\in \mathbb{C}$ depend on $\Gamma$. Then if $\wp$ is the Weierstrass elliptic function associated to $\Gamma$, the map $\mathbb{C}/\Gamma\to \mathbb{P}_2$, $\zeta\mapsto [\wp(\zeta),\wp'(\zeta),1]$, is a biholomorphism onto $C$, where $[\wp(0),\wp'(0),1]$ is interpreted as $[0,1,0]$. Since $\mathbb{C}/\Gamma$ is a complex Lie group, this biholomorphism induces a complex Lie group structure on $C$.

In this section we shall prove the following theorem.\\

\e{thm}\label{RuzzyDoo}
There exists a $6$-sheeted connected branched covering space $N$ of $\mathbb{P}_2\backslash C$ that is an Oka manifold.
\ee{thm}

Recall that a complex manifold $X$ of dimension $m$ is \emph{dominable} if for some $x\in X$, there exists a holomorphic map $\mathbb{C}^m\to X$ taking $0$ to $x$ that is a local biholomorphism at $0$. This property passes down branched covering maps. Since every Oka manifold is dominable, it follows that $\mathbb{P}_2\backslash C$ is dominable.

The branched covering space $N$ was first constructed by Buzzard and Lu \cite[Section 5.1]{Blu} as the graph complement of a meromorphic function. They showed that $N$ is dominable, implying that $\mathbb{P}_2\backslash C$ is dominable. Hanysz \cite[Theorem 4.6]{Alex} improved this result by showing that $N$ is Oka. Many of the details of the proof are omitted from both \cite{Blu} and \cite{Alex}. We will present a detailed proof of this result here based on Buzzard and Lu's original construction.

If $L$ is any line in $\mathbb{P}_2$, then Bezout's theorem implies that $L$ intersects $C$ at three points, counted with multiplicity. The group law on the cubic implies that if $p_1$, $p_2$, $p_3$ are the intersection points and
\[
p_1=[\wp(\zeta_1),\wp'(\zeta_1),1],\quad p_2= [\wp(\zeta_2),\wp'(\zeta_2),1], \quad p_3=[\wp(\zeta_3),\wp'(\zeta_3),1],
\]
then $\zeta_1+\zeta_2+\zeta_3=0$, where this addition is the group addition on $\mathbb{C}/\Gamma$. Conversely, given $\zeta_1,\zeta_2,\zeta_3\in \mathbb{C}/\Gamma$ such that $\zeta_1+\zeta_2+\zeta_3=0$, then the points $p_1$, $p_2$, $p_3$ defined as above are the intersection points of some line $L$ in $\mathbb{P}_2$ with $C$. See \cite[Section V.4]{FB} for details.

If $L$ is tangent to $C$ at $p=[x,y,z]$, then $L$ intersects $C$ at $p$ with multiplicity at least $2$. Furthermore, $L$ is the set of points $[X,Y,Z]\in\mathbb{P}_2$ such that
\e{equation}\label{Line}
F_x(x,y,z)X+F_y(x,y,z)Y+F_z(x,y,z)Z=0.
\ee{equation}

If $L$ intersects $C$ at $p$ with multiplicity $3$, then $p$ is an inflection point. The inflection points of $C$ occur where the determinant of the Hessian of $F$ vanishes. In particular, $C$ has $9$ distinct inflection points \cite[Section 4.5]{Gerdy}.

Let $E$ be the set of points $([x,y,z],[X,Y,Z])\in \mathbb{P}_2\times \mathbb{P}_2$ satisfying both (\ref{Cubic}) and (\ref{Line}). Then $(p,q)\in E$ if and only if $p\in C$ and $q$ is on the tangent line to $C$ at $p$.\\

\e{prop}
The set $E$ is a smooth, connected, $2$-dimensional, closed algebraic subvariety of $\mathbb{P}_2\times\mathbb{P}_2$.
\ee{prop}
\e{proof}
Let 
\[
A=\Set{([x,y,z],[X,Y,Z])\in \mathbb{P}_2\times \mathbb{P}_2}{z,Z\neq0}.
\]
Then $A$ is a chart on $\mathbb{P}_2\times\mathbb{P}_2$ and can be identified with $\mathbb{C}^2\times \mathbb{C}^2$ by setting $z,Z=1$. Let $H: A\to \mathbb{C}$ be defined by
\[
H(x,y,X,Y)=F_x(x,y,1)X+F_y(x,y,1)Y+F_z(x,y,1).
\]

Then $E\cap A$ is the zero set of the function $G: A\to \mathbb{C}^2$,
\[
G(x,y,X,Y)=(F(x,y,1),H(x,y,X,Y)).
\]
The Jacobian matrix of $G$ is
\[
J(G)=\left(\begin{array}{cccc}F_x & F_y & 0 & 0 \\H_x & H_y & H_X & H_Y \end{array}\right)=\left(\begin{array}{cccc}F_x & F_y & 0 & 0 \\H_x & H_y & F_x & F_y \end{array}\right).
\]
We claim that $J(G)$ has full rank at every point $([x, y, 1], [X, Y, 1])\in E\cap A$. Indeed, since $F_xX+F_yY+F_z=0$, both $F_x$ and $F_y$ cannot vanish at $[x, y, 1]$; otherwise $F_z$ would also vanish and $C$ would be singular. So we are done.
\ee{proof}

Let $\pi: E\to C, (p, q)\mapsto p$, be the projection onto $C$. Then $\pi$ is holomorphic and the fibre over $p\in C$ is the tangent line of $C$ at $p$.

Let $U_x$ be the open subset of $C$ on which $F_x$ does not vanish. Let $h_x: \pi^{-1}(U_x)\to U_x\times \mathbb{P}$, $(p,[X,Y,Z])\mapsto (p,[Y,Z])$. Then $h_x$ is a biholomorphism and hence a local trivialisation of $\pi$. The subsets $U_y, U_z$ and the corresponding biholomorphisms $h_y$, $h_z$ are defined analogously. Since $F_x$, $F_y$, $F_z$ never vanish simultaneously, we see that the sets $\pi^{-1}(U_x)$, $\pi^{-1}(U_y)$, $\pi^{-1}(U_z)$ form an open cover of $E$. Thus $E$ is a holomorphic $\mathbb{P}$-bundle over $C$.

Let $L$ be tangent to $C$ at $p$ and let $p'=-2p\in C$. Then $p'=p$ if and only if $p$ is an inflection point. Otherwise, $p'\neq p$ and $L$ intersects $C$ with multiplicity $2$ at $p$ and $1$ at $p'$.

Let $\sigma_1, \sigma_2: C\to E$, $\sigma_1 (p)= (p,p)$, $\sigma_2 (p)=(p,p')$. Then $\sigma_1, \sigma_2$ are holomorphic sections of $\pi$.

Let $u: \mathbb{C}\to C$ be the universal covering map. Then the pullback bundle $\tilde{E}$ of $E$ by $u$ is
\[
\tilde{E}= \Set{(\zeta,w)\in \mathbb{C}\times E}{u(\zeta)=\pi (w)}.
\]
Then $\tilde{E}$ is a complex manifold. Let $\tilde{u}: \tilde{E}\to E$, $(\zeta,w)\mapsto w$, and $\tilde{\pi}:\tilde{E}\to \mathbb{C}$, $(\zeta,w)\mapsto \zeta$, be the projection maps. Then $\tilde{u}$ and $\tilde{\pi}$ are holomorphic and form the following commuting square:
\[
\begin{CD}
\tilde{E} @>\tilde{u}>> E\\
@VV \tilde{\pi} V @VV \pi V\\
\mathbb{C} @> u >> C
\end{CD}
\]

Now, $u^{-1}(U_x)$, $u^{-1}(U_y)$, $u^{-1}(U_z)$ form an open cover of $\mathbb{C}$. We see that the map $\tilde{h}_x: \tilde{\pi}^{-1}(u^{-1}(U_x))\to u^{-1}(U_x)\times \mathbb{P}$, $(\zeta,w)\mapsto (\zeta,\text{proj}_2(h_x (w)))$, is a local trivialisation of $\tilde{\pi}$. Defining the maps $\tilde{h}_y$, $\tilde{h}_z$ analogously, we see that $\tilde{E}$ is a holomorphic $\mathbb{P}$-bundle over $\mathbb{C}$. The structure group of $\tilde{E}$ is the M\"obius group, which is connected. Then $\tilde{E}$ is a trivial bundle, since every holomorphic fibre bundle over a non-compact Riemann surface whose structure group is connected is holomorphically trivial \cite[Theorem 1.0]{Trivvy}.

Let $s_1,s_2 : \mathbb{C}\to \tilde{E}$, $s_1(\zeta)=(\zeta, \sigma_1(u(\zeta)))$, $s_2 (\zeta)=(\zeta, \sigma_2 (u(\zeta)))$. Then $s_1,s_2$ are well-defined holomorphic sections of $\tilde{\pi}$. The triviality of $\tilde{E}$ implies that $s_1(\mathbb{C})$, $s_2(\mathbb{C})$ may be identified with the graphs of meromorphic functions on $\mathbb{C}$.

We shall now show that $\tilde{E}\backslash s_1(\mathbb{C})$ is an affine bundle over $\mathbb{C}$. Let $h: \tilde{E}\to \mathbb{C}\times\mathbb{P}$ be a global trivialisation of $\tilde{E}$. Let $y=\text{proj}_2\circ h\circ s_1:\mathbb{C}\to\mathbb{P}$. Let $B_0=\Set{\zeta\in\mathbb{C}}{y(\zeta)\neq 0}$ and $B_\infty=\Set{\zeta\in \mathbb{C}}{y(\zeta)\neq \infty}$. Then the maps
\[
\Phi_0: (B_0\times \mathbb{P})\backslash h\circ s_1(B_0) \to B_0\times \mathbb{C}, \left(\zeta, w\right)\mapsto \left(\zeta, \frac{w}{1-\frac{w}{y(\zeta)}}\right),
\]
and
\[
\Phi_\infty: (B_\infty\times \mathbb{P})\backslash h\circ s_1(B_\infty)\to B_\infty\times \mathbb{C}, \left(\zeta, w\right)\mapsto \left(\zeta, \frac{1}{w -y(\zeta)}\right),
\]
are biholomorphisms. It follows that the maps $\Phi_0\circ h|\tilde{\pi}^{-1}(B_0)\backslash s_1(B_0)$ and $\Phi_\infty\circ h|\tilde{\pi}^{-1}(B_\infty)\backslash s_1(B_\infty)$ are local trivialisations of $\tilde{E}\backslash s_1(\mathbb{C})$. The transition maps are easily checked to be affine.

Thus $\tilde{E}\backslash s_1(\mathbb{C})$ is a $\mathbb{C}$-bundle over $\mathbb{C}$ with connected structure group, hence is a trivial bundle. We can then identify $\tilde{E}\backslash s_1(\mathbb{C})$ with $\mathbb{C}^2$. Then $\mathbb{C}^2\backslash s_2(\mathbb{C})$ is the complement of the graph of a meromorphic function.

Hanysz \cite[Theorem 4.6]{Alex} proved the following theorem.\\

\e{thm}[Hanysz]
Let $X$ be a complex manifold, and let $m:X\to\mathbb{P}$ be a holomorphic map with graph $\Gamma$. Suppose $m$ can be written in the form $m=f+1/g$ for holomorphic functions $f$ and $g$. Then $(X\times \mathbb{C})\backslash \Gamma$ is Oka if and only if $X$ is Oka.
\ee{thm}

It is well known that $\mathbb{C}$ is Oka. Furthermore, Buzzard and Lu \cite[Propositon 5.1]{Blu} use the theorems of Mittag-Leffler and Weierstrass to prove that every holomorphic map $m:\mathbb{C}\to \mathbb{P}$ can be written in the form $m=f+1/g$. Hence the manifold $\tilde{E}\backslash (s_1(\mathbb{C})\cup s_2(\mathbb{C}))$ is biholomorphic to $\mathbb{C}^2\backslash s_2(\mathbb{C})$ and satisfies the hypothesis of the theorem and so is Oka.

The map $\tilde{u}$ is the pullback of an infinite covering map along a continuous map and so is an infinite covering map $\tilde{E}\to E$. Recall that $\sigma_1(C)=\Set{(p,p)\in\mathbb{P}_2\times \mathbb{P}_2}{p\in C}$, $\sigma_2(C)=\Set{(p,p')\in\mathbb{P}_2\times \mathbb{P}_2}{p\in C}$. By definition of the pullback bundle, the preimage $\tilde{u}^{-1}(\sigma_1(C)\cup \sigma_2(C))$ consists of pairs $(\zeta,(p,p))$ or $(\zeta,(p,p'))$, where $u(\zeta)=p$. But by definition of $s_1, s_2$, if $u(\zeta)=p$, then $s_1(\zeta)= (\zeta,(p,p))$ and $s_2(\zeta)=(\zeta, (p,p'))$, implying $\tilde{u}^{-1}(\sigma_1(C))=s_1(\mathbb{C})$ and $\tilde{u}^{-1}(\sigma_2(C))=s_2(\mathbb{C})$. Hence the restriction $\tilde{u}|\tilde{E}\backslash(s_1(\mathbb{C})\cup s_2(\mathbb{C}))$ is an infinite covering map $\tilde{E}\backslash (s_1(\mathbb{C})\cup s_2(\mathbb{C})) \to E\backslash (\sigma_1 (C)\cup \sigma_2(C))$. It follows that $N=E\backslash (\sigma_1 (C)\cup \sigma_2(C))$ is Oka \cite[Proposition 5.5.2]{Forstneric}. Since $E$ is connected, so is $N$.

Let $\lambda:E\to\mathbb{P}_2$, $(p,q)\mapsto q$. Then $\lambda$ is holomorphic and, since $E$ is compact, is proper. The preimage $\lambda^{-1}(C)$ consists precisely of those $(p,q)\in \mathbb{P}_2\times\mathbb{P}_2$ such that $q$ is in the intersection of the tangent line to $C$ at $p$ with $C$ itself. Hence either $q=p$ or $q=p'$. Therefore $\lambda^{-1}(C)=\sigma_1(C)\cup \sigma_2(C)$ and the restriction $\lambda|{N}$ is a proper holomorphic map from an Oka manifold into $\mathbb{P}_2\backslash C$.

We shall now show that $\lambda|{N}\to\mathbb{P}_2\backslash C$ is a surjective finite map. Given $[q_1,q_2,q_3]\in\mathbb{P}_2\backslash C$, the curve defined by
\[
q_1F_x+q_2F_y+q_3F_z=0
\]
is of degree $2$ and by Bezout's theorem intersects $C$ at $6$ points, counted with multiplicity. It is clear that these $6$ intersection points are precisely the points on $C$ whose tangents pass through $q$. Then $\lambda$ is surjective and finite. Furthermore, the fibres generically contain $6$ points \cite[Section 5.7]{Gerdy} and so $\lambda$ is a $6$-sheeted branched covering of $\mathbb{P}_2\backslash C$ by an Oka manifold. Note that $\lambda$ does have branching as the curve
\[
q_1F_x+q_2F_y+q_3F_x=0
\]
intersects $C$ with multiplicity $2$ precisely at the inflection points of $C$. It follows that the set of critical values of $\lambda$ is the union of the $9$ inflection tangents of $C$ with the inflection points removed. Indeed Bezout's theorem implies that each inflection tangent intersects $C$ only at its corresponding inflection point. Then the intersection of any two inflection tangents occurs in $\mathbb{P}_2\backslash C$ and so the set of critical values forms a singular hypersurface of $\mathbb{P}_2\backslash C$ whose singular locus consists of the intersection points of the inflection tangents.


\section{The Main Theorem}\label{MThm}

Let $X$ be a compact Riemann surface of genus $1$ and $R_3$ the space of holomorphic maps $X\to\mathbb{P}$ of degree $3$. Recall that $R_3$ is a connected $6$-dimensional complex manifold. We can now state our main result.\\

\e{thm}\label{main}
There exists a $6$-sheeted connected branched covering space of $R_3$ that is an Oka manifold.
\ee{thm}

To prove this theorem, we will make use of the following construction. Recall that the map $s_3: X\to X$, $x\mapsto 3x$, is a $9$-sheeted unbranched covering map. Let $\chi$ be the map $R_3/M\to X$ induced by the map $\beta: R_3\to X$, $f\mapsto \psi(D_0(f))$, where $\psi$ is the Jacobi map $S^3X\to X$, $p_1+p_2+p_3\mapsto p_1+p_2+p_3$, and $D_0: R_3\to S^3X$ takes $f$ to its divisor of zeroes. Proposition \ref{Top Cover} implies that, for a certain smooth cubic $C$ in $\mathbb{P}_2$ biholomorphic to $X$ (called $C_0$ before), there exists an unbranched $9$-sheeted covering map $\tilde{s}_3: \mathbb{P}_2\backslash C\times X\to R_3/M$ such that the following diagram is a pullback square.
\[
\begin{tikzcd}
\mathbb{P}_2\backslash C \times X \arrow{r}{\tilde{s}_3} \arrow{d}{\text{proj}_2} & R_3/M \arrow{d}{\chi} \\
X \arrow{r}{s_3} & X
\end{tikzcd}
\]
Furthermore, Proposition \ref{Cubic Complement} allows us to identify $\mathbb{P}_2\backslash C$ with $\chi^{-1}(0)$. Let $\pi: R_3\to R_3/M$ be the quotient map. Then $\pi(\beta^{-1}(0))=\chi^{-1}(0)$. Recall that if $f\in \beta^{-1}(0)$, then $\tilde{s}_3(\pi(f),x)=\pi(f\circ\tau_{-x})$. For the remainder of this section we shall use this identification to denote all elements in $\mathbb{P}_2\backslash C$ by $\pi(f)$, where the representative $f\in \beta^{-1}(0)$ is chosen appropriately. As in Section \ref{OBC}, let
\[
N=\Set{(p,q)\in C\times\mathbb{P}_2\backslash C}{\text{$q$ is in the tangent to $C$ at $p$}},
\]
and take $\lambda: N\to \mathbb{P}_2\backslash C$, $(p,q)\mapsto q$. Then by Theorem \ref{RuzzyDoo}, $N$ is an Oka manifold and $\lambda$ is a $6$-sheeted branched covering map.

Let $S$ be the $9$-element subgroup $\Set{x\in X}{3x=0}$, isomorphic to $\mathbb{Z}_3\times \mathbb{Z}_3$. (Note that if $\eta_0:X\to C$ is the embedding defining $C=C_0$ as in Section \ref{Rum}, then $\eta_0(S)$ consists of the inflection points of $C$.) Then $S$ acts on $\mathbb{P}_2\backslash C\times X$  via the action 
\[
S\times\mathbb{P}_2\backslash C\times X\to \mathbb{P}_2\backslash C\times X, \quad t\cdot(\pi(f),x)=(\pi(f\circ\tau_t),x+t).
\]
It is clear that this action realises $S$ as the group of covering transformations of $\tilde{s}_3$.

We will now lift the above action of $S$ on $\mathbb{P}_2\backslash C\times X$ to an action on $N\times X$. Recall that $C$ is the image of an embedding $p$ of $X$ into $\mathbb{P}_2$. We claim that the map
\[
S\times N\times X\to N\times X, \quad t \cdot ((p(y),\pi(f)), x)= ((p(y-t), \pi(f\circ\tau_t)), x+t),
\]
is an action of $S$ on $N\times X$ such that whenever $t\in S$,
\[
(\lambda\times\id_X)(t\cdot((p,q),x))=t\cdot((\lambda\times\id_X)((p,q),x)).
\]
To check this, it is enough to know that whenever $t\in S$ and $(p(y), \pi(f))\in N$, $t\cdot(p(y), \pi(f))\in N$. To show that $t\cdot(p(y), \pi(f))=(p(y-t), \pi(f\circ\tau_{t}))\in N$, we must firstly show that $f\circ\tau_{t}\in\beta^{-1}(0)$. This follows from the fact that $\beta(f\circ\tau_t)=\beta(f)-3t=0-0=0$. It remains to show that $\pi(f\circ\tau_t)$ is on the tangent to $C$ at $p(y-t)$. We know that $\pi(f)$ is on the tangent to $C$ at $p(y)$. Then by Proposition \ref{LineLemma}, $y+y+(-2y)$ is a divisor corresponding to a fibre of $f$. But then $(y-t)+(y-t)+(-2(y-t))$ is a divisor corresponding to a fibre of $f\circ\tau_t$ and so $\pi(f\circ\tau_t)$ is on the tangent to $C$ at $p(y-t)$. Hence we have defined an action of $S$ on $N\times X$.

Note that it is generally not possible to lift a group action up a branched covering map. It works in our special case because $N\subset C\times \mathbb{P}_2\backslash C$, so lifting an action from $N\times X$ to $C\times \mathbb{P}_2\backslash C\times X$ amounts to choosing an action on $C$ that restricts to an action on $N\times X$. Indeed we have the following commuting diagram.
\[
\begin{tikzcd}
N\times X \arrow[hookrightarrow]{r} \arrow{d}{\lambda\times \id_X}& C\times (\mathbb{P}_2\backslash C)\times X \arrow{dl}{\proj}\\
\mathbb{P}_2\backslash C\times X
\end{tikzcd}
\]

We will now use the action of $S$ on $N\times X$ to construct a $6$-sheeted Oka branched covering space of $R_3/M$. Indeed, the action of $S$ on $X$ is free so the action of $S$ on $N\times X$ is free. It follows that the quotient $(N\times X)/S$ is a complex manifold such that the quotient map $\zeta:N\times X\to (N\times X)/S$ is a holomorphic $9$-sheeted unbranched covering. Since $N$ and $X$ are Oka, so is $N\times X$. Then $(N\times X)/S$ is also Oka.

Let $\Psi=\tilde{s}_3\circ(\lambda\times \text{id}_X): N\times X\to R_3/M$. Then $\Psi$ is the composition of a $6$-sheeted branched covering followed by a $9$-sheeted unbranched covering map. We now claim that $\Psi$ is $S$-invariant. Indeed,
\begin{align*}
\Psi(t\cdot((p(y),\pi(f)),x))&=\tilde{s}_3((\lambda\times \text{id}_X)(t\cdot((p(y),\pi(f)),x)))\\
&= \tilde{s}_3(t\cdot((\lambda\times \text{id}_X)((p(y),\pi(f)),x)))\\
&= \tilde{s}_3((\lambda\times \text{id}_X)((p(y),\pi(f)),x)).
\end{align*}
Thus $\Psi$ induces a holomorphic map $\Lambda:(N\times X)/S\to R_3/M$ such that the following diagram commutes.
\[
\begin{tikzcd}
N\times X \arrow{dr}{\Psi}\arrow{r}{\zeta} \arrow{d}{\lambda\times\id_X}& (N\times X)/S \arrow{d}{\Lambda}\\
\mathbb{P}_2\backslash C\times X \arrow{r}{\tilde{s}_3} & R_3/M
\end{tikzcd}
\]
Since $\Psi$ is a finite map, $\Lambda$ is also a finite map. Moreover, it is clear that $S$ acts freely on the fibres of $\Psi$, which generically contain $54$ points. It follows that the fibres of $\Lambda$ generically contain $6$ points and so $\Lambda$ is a $6$-sheeted branched covering map.

To prove Theorem \ref{main}, let $W$ be the pullback of $R_3$ along $\Lambda$. Then we obtain the following pullback square:
\[
\begin{tikzcd}
W \arrow{r}{\Phi} \arrow{d}{\tilde{\pi}}& R_3 \arrow{d}{\pi}\\
(N\times X)/S \arrow{r}{\Lambda} & R_3/M
\end{tikzcd}
\]
where $\tilde{\pi}: W \to (N\times X)/S$ is the pullback of $\pi$ along $\Lambda$ and $\Phi: W\to R_3$ is the pullback of $\Lambda$ along $\pi$. Then $W$ is smooth since $\pi$ is a submersion and is a $6$-sheeted branched covering space of $R_3$. By Theorem \ref{AutY}, $\pi$ is a principal $M$-bundle, which implies that $\tilde{\pi}$ is also a principal $M$-bundle. Since $M, N$ and $X$ are connected, so is $W$. Now, since $M$ is a complex Lie group, $M$ is Oka. Then since $(N\times M)/S$ is Oka it follows that $W$ is Oka and the theorem is proved.

Note that $\Phi$ has branching because $\lambda$ has branching, as mentioned at the end of Section \ref{OBC}. The set of critical points of $\Phi$ will be described in Section \ref{Walt Disney}.

Recall that a complex manifold of dimension $m$ is \emph{dominable} if for some $x\in X$, there is a holomorphic map $\mathbb{C}^m\to X$ taking $0$ to $x$ that is a local biholomorphism at $0$. We say $X$ is \emph{$\mathbb{C}$-connected} if given any two points $x,y\in X$ there exists a holomorphic map $\mathbb{C}\to X$ whose image contains both $x$ and $y$. It is clear that both properties pass down branched covering maps. Every Oka manifold is known to be both dominable and $\mathbb{C}$-connected. It follows that $R_3$ is both dominable and $\mathbb{C}$-connected.\\

\e{cor}
The complex manifold $R_3$ is dominable and $\mathbb{C}$-connected.
\ee{cor}

\chapter{Further Results\label{ch:three}}

\section{Strong Dominability}\label{StrongDom}

Recall that a complex manifold $Y$ of dimension $n$ is \emph{strongly dominable} if for every $y\in Y$, there is holomorphic map $f:\mathbb{C}^{n}\to Y$ such that $f(0)=y$ and $f$ is a local biholomorphism at $0$.

Let $X$ be a compact Riemann surface of genus $1$ and $R_3$ the corresponding complex manifold of holomorphic maps $X\to\mathbb{P}$ of degree $3$. Recall that $X$ is biholomorphic to a torus $\mathbb{C}/\Gamma$, where $\Gamma$ is a lattice in $\mathbb{C}$. Let $\Gamma_0$ denote the hexagonal lattice. Then Proposition \ref{hex} implies that $\mathbb{C}/\Gamma_0$ is the unique torus whose group of automorphisms fixing the identity has order $6$. The goal of this section is to show that if $X$ is not biholomorphic to $\mathbb{C}/\Gamma_0$, then the complex manifold $R_3$ is strongly dominable.

Let $C$ be a smooth cubic curve in $\mathbb{P}_2$ biholomorphic to $X$ and defined by the equation $F(x,y,z)=0$. Let $N$ be the set of $(p,q)\in\mathbb{P}_2\times(\mathbb{P}_2\backslash C)$ such that $p\in C$ and $q$ is on the line tangent to $C$ at $p$. We know from Section \ref{OBC} that $N$ is an Oka manifold. Furthermore, the map $\lambda: N\to \mathbb{P}_2\backslash C$, $(p,q)\mapsto q$, is a $6$-sheeted branched covering map.

Suppose that $X$ is not biholomorphic to $\mathbb{C}/\Gamma_0$. We will show that $R_3$ is strongly dominable by showing that $\mathbb{P}_2\backslash C$ is strongly dominable. To prove the latter, it is sufficient to know that every fibre of $\lambda$ contains a regular point. To see why this implies that $R_3$ is strongly dominable, let $S$ be the set of $x\in X$ such that $3x=0$ and $M$ be the M\"obius group. Recall from the proof of Theorem \ref{main} the maps $\Psi: N\times X\to R_3/M$ and $\Lambda: (N\times X)/S\to R_3/M$. Recall that $\Psi$ is $\lambda\times \text{id}_X$ postcomposed by a $9$-sheeted unbranched covering map. Every fibre of $\lambda$ contains a regular point if and only if every fibre of $\Psi$ contains a regular point. The map $\Lambda$ is a $6$-sheeted branched covering map from Section \ref{MThm}, while the quotient map $\zeta: N\times S\to (N\times S)/S$ is a $9$-sheeted unbranched covering map such that $\Lambda\circ \zeta=\Phi$. Hence every fibre of $\Lambda$ contains a regular point if and only if every fibre of $\Psi$ contains a regular point. Let $W$ and $\Phi:W\to R_3$ be as in the proof of Theorem \ref{main}. Then $\Phi$ is the pullback of $\Lambda$ along the quotient map $\pi: R_3\to R_3/M$, so every fibre of $\Phi$ contains a regular point if and only if every fibre of $\Lambda$ contains a regular point. Hence both $R_3$ and $\mathbb{P}_2\backslash C$ will be strongly dominable if we can show every fibre of $\lambda$ contains a regular point.\\

\e{lem}
A point $(p,q)\in N$ is a critical point of $\lambda$ if and only if $p$ is an inflection point of $C$.
\ee{lem}
\e{proof}
Fix $(p,q)\in N$ and let $[x,y,z]$ denote homogeneous coordinates on $\mathbb{P}_2$. Using a projective transformation we can assume that $p=[0,0,1]$ and $q=[1,0,0]$ and so the tangent line to $C$ at $p$ is the line $y=0$. Then $F_x(p)=F_z(p)=0$ and $F_{y}(p)\neq0$. Let $A\subset \mathbb{P}_2$ be the affine chart obtained by setting $z=1$. It follows by the implicit function theorem that there exists a neighbourhood $\Omega$ of $p$ in $A$, a neighbourhood $U$ of $0$ in $\mathbb{C}$ and a holomorphic map $f:U\to \mathbb{C}$ such that $C\cap \Omega$ is the graph of $f$ via the embedding $U\to \Omega$, $t\mapsto (t,f(t))$. Then $f(0)=0$ and, since $y=0$ is tangent to $C$ at $p$, $f'(0)=0$.

We claim that $(p,q)$ is a critical point of $\lambda$ if and only if $f''(0)=0$. Indeed, if $f''(0)=0$, then $f'$ is not injective on any neighbourhood of $p$. It follows that for every neighbourhood $V$ of $p$ in $\Omega$, there exist two distinct points in $V\cap C$ whose tangent lines are parallel in the affine chart $A$ and hence intersect on the line at infinity. Then given any neighbourhood $V'$ of $q$ in $\mathbb{P}_2$ we can choose $V$ sufficiently small that the tangent line of every point in $V$ intersects the line at infinity in $V'$. Then we can find two distinct points in $V$ whose tangent lines intersect in $V'$. Hence given any neighbourhood of the form $V\times V'$ of $(p,q)$ in $N$, we can find two distinct points in $V\times V'$ in the same fibre of $\lambda$ and so $(p,q)$ is a critical point of $\lambda$.

Conversely, suppose that $f''(0)\neq0$. Then shrinking $U$ and $\Omega$ if necessary, we can assume that $f'$ is injective on $U$. It is easy to check that for distinct points $t_1, t_2\in U$, the tangent lines to $t_1, t_2$ intersect in the affine chart $A$ at the point
\[
(X,Y)=\left(\frac{f'(t_1)t_1-f'(t_2)t_2+f(t_2)-f(t_1)}{f'(t_1)-f'(t_2)}, f'(t_1)X+f(t_1)-f'(t_1)t_1\right).
\]
To show that $(p,q)$ is a regular point of $\lambda$, it is sufficient to find a neighbourhood on which $\lambda$ is injective. To show this, it is enough to know that given $\epsilon>0$ we can find $\delta>0$ such that  $|X|<\epsilon$ whenever $|t_1|<\delta$ and $|t_2|<\delta$. Indeed, since $Y$ depends continuously on $X$, this will show that given two sequences $(r_n)$, $(s_n)$ of points on $C$ converging to $p$ such that $r_n\neq s_n$ for all $n$, the sequence of intersection points $(X_n, Y_n)$ will also converge to $p$. Take $\epsilon>0$. We see that
\begin{align*}
|X|&=\left|\frac{f'(t_1)t_1-f'(t_2)t_1+f'(t_2)t_1-f'(t_2)t_2+f(t_2)-f(t_1)}{f'(t_1)-f'(t_2)}\right|\\
&\leq |t_1|+|f'(t_2)|\left|\frac{t_1-t_2}{f'(t_1)-f'(t_2)}\right|+\left|\frac{f(t_2)-f(t_1)}{f'(t_1)-f'(t_2)}\right|.
\end{align*}
Restrict $t_1$ such that $|t_1|<\epsilon$. Since $f'(0)=0$, for $|t_2|$ sufficiently small, $|f'(t_2)|<\epsilon$. Furthermore, since $f'$ is injective on $U$, $f''(0)\neq0$ and, restricting $t_1$, $t_2$ further if necessary, we see that $|t_1-t_2|/|f'(t_1)-f'(t_2)|$ is bounded above by some $M>0$. Finally,
\[
\left|\frac{f(t_2)-f(t_1)}{f'(t_1)-f'(t_2)}\right|=\left|\frac{t_1-t_2}{f'(t_1)-f'(t_2)}\right|\cdot\left|\frac{f(t_1)-f(t_2)}{t_1-t_2}\right|<M\epsilon,
\]
where $t_1$, $t_2$ are restricted so that $|f(t_1)-f(t_2)|/|t_1-t_2|<\epsilon$. Then we see that
\[
|X|<\epsilon+2M\epsilon,
\]
which implies that $\lambda$ is injective on a neighbourhood of $(p,q)$, so $(p,q)$ is a regular point.

To complete the proof of the lemma, we must show that $p$ is an inflection point of $C$ if and only if $f''(0)=0$. Observe that for all $t\in U$, $F(t,f(t),1)=0$. Differentiating twice, we see that
\begin{multline*}
F''(t,f(t),1)=F_{xx}(t,f(t),1)+F_{xy}(t,f(t),1)f'(t)\\+(F_{xy}(t,f(t),1)+F_{yy}(t,f(t),1)f'(t))f'(t)+F_y(t,f(t),1)f''(t)=0.
\end{multline*}
Since $f'(0)=0$, evaluating at $t=0$ yields $f''(0)=-F_{xx}(p)/F_{y}(p)$ and we are done if we can show that $p$ is an inflection point if and only if $F_{xx}(p)=0$.

Let $G(x)=F(x,0,1)$. Then $G$ has the following Taylor expansion about the point $p$:
\[
G(x)=F(p)+x F_x(p)+x^2F_{xx}(p)+\cdots.
\]
Since $C$ is of degree $3$, the point $p$ is an inflection point of $C$ if and only if $G(x)$ has a zero of order $3$ at $0$ \cite[Section 2.5]{Gerdy}. Since $p\in C$, $F(p)=0$. Furthermore, since the line $y=0$ is tangent to $C$, $F_x(p)=0$. Thus $p$ is an inflection point if and only if $F_{xx}(p)=0$.
\ee{proof}

It follows from the above lemma that if $q\in \mathbb{P}_2\backslash C$, then $\lambda^{-1}(q)$ contains no regular point only if there exist three inflectional tangents to $C$ that meet at $q$. The existence of three inflectional tangents with a common point depends on the cubic $C$. We have been unable to find a reference that explains for which cubics this occurs, so we will summarise our own work.                                                       

The following lemma is well known \cite[Lemma 2.1]{Hesse}.\\

\e{lem}
Every smooth cubic in $\mathbb{P}_2$ is projectively equivalent to a cubic defined by the equation $x^3+y^3+z^3+t xyz=0$, for some $t\in\mathbb{C}$.
\ee{lem}

Let $C_t$ denote the cubic defined by the equation $x^3+y^3+z^3+t xyz=0$. Precomposing $\lambda$ by a projective transformation, we can assume that $C=C_t$ for some $t\in\mathbb{C}$. Let $\varepsilon=e^{2i\pi/3}$. The inflection points of $C_t$ are independent of $t$ and are as follows \cite[p. 238]{Hesse}:
\[
\begin{array}{ccc}
[0, 1, -1] &[0,1,-\varepsilon] &[0,1,-\varepsilon^2]\\\relax
[1,0,-1] &[1,0,-\varepsilon^2] &[1,0,-\varepsilon]\\\relax
[1,-1,0] &[1,-\varepsilon, 0] &[1,-\varepsilon^2,0].
\end{array}
\]

The tangents to $C_t$ at the above inflection points have dual coordinates as follows:
\[
\begin{array}{ccc}
[-t,3,3] &[3, -t, 3] & [3, 3, -t]\\\relax
[-t\varepsilon, 3, 3\varepsilon^2] & [3, -t\varepsilon^2, 3\varepsilon] &[3. 3\varepsilon^2, -t\varepsilon]\\\relax
[-t\varepsilon^2, 3, 3\varepsilon] &[3, -t\varepsilon, 3\varepsilon^2] &[3, 3\varepsilon, -t\varepsilon^2]
\end{array}
\]

Three inflectional tangents intersect in a point if and only if the corresponding dual points are collinear. Considering each of the $84$ cases, the only values of $t$ for which this occurs are the following:

\[
\begin{array}{ccc}
\begin{aligned}
t&=0\\
&\empty \\
&\empty
\end{aligned} &
\begin{aligned}
t &= 6\\
t&=6\varepsilon\\
 t&=6\varepsilon^2
\end{aligned} &
\begin{aligned}
t &= -3\\
t &= -3\varepsilon\\
t &= -3\varepsilon^2
\end{aligned} 
\end{array}
\]

The cubics corresponding to $t=-3, -3\varepsilon, -3\varepsilon^2$ are singular. These are the only values of for which $C_t$ is singular \cite[p. 241]{Hesse}.

Two cubics in $\mathbb{P}_2$ are projectively equivalent if and only if they are biholomorphic \cite[p.~127]{Dolgo}. Two values of $t$ give the same cubic up to biholomorphism if and only if they are in the same fibre of the map
\[
j:\mathbb{C}\to \mathbb{P}, \quad t\mapsto [2^3t^3(1-t^3/6^3)^3, 3^3(1+t^3/3^3)^3]
\]
\cite[Equation 3.11, p.~133]{Dolgo}. This map is the well-known $j$-invariant. It is evident that $j(t)=j(\varepsilon t)$. Futhermore, it is clear that that $j(6)=j(0)$ and so up to biholomorphism, there is only one smooth cubic with three inflectional tangents sharing a point. This cubic, $C_0$, is biholomorphic to $\mathbb{C}/\Gamma_0$, where $\Gamma_0$ is the hexagonal lattice, and is called the \emph{equianharmonic cubic} \cite[Theorem 3.1.3, Definition 3.1.2]{Dolgo}.

We have proved the following theorem.\\

\e{thm}
Let $C$ be a smooth cubic curve in $\mathbb{P}_2$ that is not the equianharmonic cubic. Then $\mathbb{P}_2\backslash C$ is strongly dominable.

Let $X$ be a compact Riemann surface of genus $1$ and $R_3$ the complex manifold of degree $3$ holomorphic maps $X\to \mathbb{P}$. Then if $X$ is not biholomorphic to the equianharmonic cubic, $R_3$ is strongly dominable.
\ee{thm}

Note that if $C$ is biholomorphic to the equianharmonic cubic, it is an open problem whether $\mathbb{P}_2\backslash C$ is strongly dominable.

\section{The Remmert Reduction of $R_3$}\label{UncleRemmy}

Let $X$ be a compact Riemann surface of genus $1$ and $R_3$ the complex manifold of degree $3$ holomorphic maps $X\to \mathbb{P}$. Identify $X$ with $\mathbb{C}/\Gamma$, for some lattice $\Gamma$ in $\mathbb{C}$, making $X$ a complex Lie group. The main goal of this section is to show that $R_3$ is holomorphically convex and to determine its Remmert reduction.

If $Y$ is a holomorphically convex manifold, and we define an equivalence relation on $Y$ by setting $x\sim y$ whenever $f(x)=f(y)$ for all $f\in \mathcal{O}(Y)$, then $Y/\!\!\sim$ is a Stein space such that every holomorphic map from $Y$ to a holomorphically separable space $Z$ admits a unique factorisation through the quotient map $\Pi:Y\to Y/\!\!\sim$ by a holomorphic map $Y/\!\!\sim\,\to Z$ \cite[Theorem 57.11]{KnK}. The space $Y/\!\!\sim$ is called the \emph{Remmert reduction} of $Y$. The map $\Pi$ is proper with connected fibres and induces an isomorphism $\mathcal{O}(Y/\!\!\sim)\to\mathcal{O}(Y)$. 

Recall that $X$ acts on $R_3$ by precomposition of translations. That is, if $t\in X$ and $\tau_t:X\to X$, $x\mapsto x+t$, then if $f\in R_3$, we define $t\cdot f=f\circ \tau_{-t}$. This action is holomorphic by Theorem \ref{Douady}, as the map $X\times R_3\times X\to \mathbb{P}$, $(t,f,x)\mapsto f(x-t)$, is obviously holomorphic in $x$ and $t$, and is holomorphic in $f$ because the evaluation map $e:R_3\times X\to \mathbb{P}$ is holomorphic. It is obvious that $R_3$ is not holomorphically separable because it contains the compact $X$-orbits.

We now claim that the action of $X$ on $R_3$ is free. Indeed, suppose that $f\in R_3$ and $t\in X$ such that $f=f\circ\tau_{-t}$. If $b\in X$ is a branch point of $f$, then $b+t$ is a branch point of $f\circ\tau_{-t}$. But then $f$ also has a branch point at $b+t$ and $f(b+t)=f(b)$. So either $t=0$ or $f$ attains the value $f(b)$ with multiplicity at least four, which is impossible if $f$ is of degree $3$. Finally, since $X$ is compact, it acts properly on $R_3$. Thus applying Theorem \ref{Holmann} we obtain the following result.\\

\e{thm}
The orbit space $R_3/X$ is a complex manifold such that the quotient map $R_3\to R_3/X$ is a principal $X$-bundle.
\ee{thm}

We now claim that the quotient $R_3/X$ is Stein. Recall from Section \ref{Rum} that $\beta: R_3\to X$ is the map taking $f\in R_3$ to the sum of its zeroes in the Lie group $X$, counted with multiplicity. Let $\pi:R_3\to R_3/M$ be the quotient map, where $M$ is the M\"obius group. Theorem \ref{AutY} implies that $\pi$ is a principal $M$-bundle. It is clear that $M$ acts on each fibre of $\beta$. Let $0$ be the identity element in $X$. Then $\pi|\beta^{-1}(0)\to \beta^{-1}(0)/M$ is also a principal $M$-bundle. By Proposition \ref{Cubic Complement}, $\beta^{-1}(0)/M$ is biholomorphic to $\mathbb{P}_2\backslash C$ for some smooth cubic curve $C$ biholomorphic to $X$, so $\beta^{-1}(0)/M$ is Stein. Since $M$ is Stein, $\pi|\beta^{-1}(0)$ is a principal bundle with Stein fibre and Stein base, and so $\beta^{-1}(0)$ is a $5$-dimensional Stein manifold \cite[Th\'eor\`eme 4]{Prince}.

Now, let $S$ be the set containing $0$ and all elements in $X$ of order $3$. Then $S$ is a $9$-element subgroup of $X$ and $S$ acts on $\beta^{-1}(0)$. Since the action of $X$ on $R_3$ is free, it follows that $S$ acts freely (and properly discontinuously) on $\beta^{-1}(0)$ and so the quotient map $\Pi':\beta^{-1}(0)\to\beta^{-1}(0)/S$ is a holomorphic $9$-sheeted covering map. Since $\beta^{-1}(0)$ is a Stein manifold, so is $\beta^{-1}(0)/S$.

Let $\iota:\beta^{-1}(0)\hookrightarrow R_3$ be the inclusion map and $\Pi:R_3\to R_3/X$ the quotient map. Then each fibre of $\Pi \circ \iota$ is a single $S$-orbit in $\beta^{-1}(0)$. Indeed, if $[\tilde{f}]\in R_3/X$, then there is $f\in \beta^{-1}(0)$ with $(\Pi\circ \iota)(f)=[\tilde{f}]$. If $g\in \beta^{-1}(0)$ such that $(\Pi\circ\iota)(g)=[\tilde{f}]$, then for some $t\in X$, $g=f\circ\tau_{-t}$. But $\beta(g\circ\tau_{-t})=3t=0$, which implies that $g$ is in the $S$-orbit of $f$. Thus $\Pi\circ \iota$ is an $S$-invariant holomorphic map $\beta^{-1}(0)\to R_3/X$ and so, by the universal property of the quotient, induces a bijective holomorphic map $\varphi: \beta^{-1}(0)/S\to R_3/X$ such that the following diagram commutes:
\[
\begin{tikzcd}
\beta^{-1}(0)\arrow{r}{\iota} \arrow{d}{\Pi'} & R_3 \arrow{d}{\Pi} \\
\beta^{-1}(0)/S \arrow{r}{\varphi} & R_3/X
\end{tikzcd}
\]
Thus, $R_3/X$ is Stein, and we have proved the following theorem.\\

\e{thm}
The complex manifold $R_3$ is holomorphically convex. Furthermore, its Remmert reduction is $R_3/X$.
\ee{thm}

\section{An Alternative Proof of the Main Theorem}\label{Walt Disney}

Let $X$ be a compact Riemann surface of genus $1$, and $R_3$ be the complex manifold of holomorphic maps $X\to\mathbb{P}$ of degree $3$. Let $M$ be the M\"obius group and identify $X$ with $\mathbb{C}/\Gamma$ for $\Gamma$ a lattice in $\mathbb{C}$, making $X$ a Lie group. Then $M$ and $X$ act freely and properly on $R_3$ and so by Theorem \ref{Holmann} the quotient maps $\pi: R_3\to R_3/M$ and $\Pi: R_3\to R_3/X$ are principal bundles with fibres $M$ and $X$ respectively. We know from Proposition \ref{Top Cover} that there is a cubic curve $C$ in $\mathbb{P}_2$ biholomorphic to $X$ and a $9$-sheeted unbranched covering map $\tilde{s}_3:\mathbb{P}_2\backslash C\times X\to R_3/M$. Theorem \ref{RuzzyDoo} states that there is a $6$-sheeted branched covering map $\lambda:N\to\mathbb{P}_2\backslash C$ such that $N$ is an Oka manifold. In Section \ref{MThm} we constructed an Oka manifold $W$ and a $6$-sheeted branched covering map $\Phi:W\to R_3$.

In this section we shall construct an alternative $6$-sheeted branched covering space of $R_3$, before showing that it is isomorphic to $W$ in a natural way. We will then determine the critical values of both $6$-sheeted branched covering maps.

Let $\beta:R_3\to X$ be the map taking $f\in {R}_3$ to the Lie group sum of its zeroes, counted with multiplicity. We showed in Section \ref{UncleRemmy} that the quotient map $\varphi:\beta^{-1}(0)\to R_3/X$ is a $9$-sheeted unbranched covering map. We know from Proposition \ref{Cubic Complement} that $\beta^{-1}(0)/M$ is biholomorphic to $\mathbb{P}_2\backslash C$. The quotient map $\beta^{-1}(0)\to \beta^{-1}(0)/M$ is the map $\pi'=\pi|\beta^{-1}(0)$ and is a principal $M$-bundle. Let $Y$ be the pullback of $\beta^{-1}(0)$ along the map $\lambda$. Then we have the following commuting diagram:
\[
\begin{tikzcd}
Y\arrow{r}{\lambda'} \arrow{d} &   \beta^{-1}(0)\arrow{d}{\pi'}\\
N \arrow{r}{\lambda} & \mathbb{P}_2\backslash C
\end{tikzcd}
\]
where $\lambda'$ is the pullback of $\lambda$ along $\pi'$. Then $Y$ is a principal $M$-bundle over $N$. Since $M$ and $N$ are Oka, so is $Y$.

Let $S\subset X$ be the $9$-element subgroup $\Set{x\in X}{3x=0}$. We know from Section \ref{MThm} that $S$ acts freely on $N\times X$, and that $W$ is an $M$-bundle over $(N\times X)/S$. It can be easily shown that this action defines an action on $N$ by projecting onto the first component. Also, $S$ acts on $Y$. Namely, for $(n,f)\in Y\subset N\times \beta^{-1}(0)$, let $t\cdot (n,f)= (t\cdot n, f\circ\tau_{t})$. Since $t\in S$ it follows that $f\circ\tau_{t}\in \beta^{-1}(0)$. To show that this defines an action of $S$ on $Y$ it is enough to show that $t\cdot (n,f)\in Y$. This is equivalent to showing that $\lambda(t\cdot n)= \pi(f\circ \tau_{t})$, which is clear from the definition of the action of $S$ on $N\times X$ in Section \ref{MThm}. Hence $S$ acts on $Y$ and because the action of $S$ on $\beta^{-1}(0)$ is free (as the action of $X$ on $R_3$ is free), we see that the quotient map $Y\to Y/S$ is a $9$-sheeted unbranched covering map and so $Y/S$ is an Oka manifold. Now $\varphi\circ \lambda'$ is a $54$-sheeted branched covering map $Y\to R_3/X$. Furthermore, $\varphi\circ\lambda'$ is $S$-invariant as $\varphi(f\circ\tau_t)=\varphi(f)$ for all $t\in S$ and $f\in \beta^{-1}(0)$. It follows that $\varphi\circ \lambda'$ induces a finite holomorphic map $\Lambda': Y/S\to R_3/X$. It is easy to see that $S$ acts freely on the fibres of $\varphi\circ\lambda'$. Hence $\Lambda'$ is a $6$-sheeted branched covering map.

Let $W'$ be the pullback of $Y/S$ along the quotient map $\Pi$ and $\Phi':W\to R_3$ be the pullback of $\Lambda'$ along $\Pi$. Then we have the following commuting diagram:
\[
\begin{tikzcd}
W'\arrow{r}{\Phi'}\arrow{d} & R_3\arrow{d}{\Pi}  \\
Y/S \arrow{r}{\Lambda'} & R_3/X
\end{tikzcd}
\]

and $\Phi'$ is $6$-sheeted branched covering map $W'\to R_3$. Furthermore, $W'$ is a principal $X$-bundle over the Oka manifold $Y/S$, so is Oka.

The following commuting diagram of complex manifolds and holomorphic maps shows the relationships between $W$ and $W'$. Branched and unbranched covering maps are labelled with the number of sheets. The unbranched coverings are precisely the maps with $9$ sheets. Furthermore we label principal $M$- and $X$-bundles with their corresponding fibre.

\[
\begin{tikzcd}
N\times X \arrow{dr}{\lambda\times \text{id}_X}[swap]{6}\arrow{dd}{X} \arrow{rr}[swap]{9}&&(N\times X)/S \arrow{d}{6}[swap]{\Lambda}& W \arrow{l}{M}\arrow{ddl}{6}[swap]{\Phi}\\
 &\mathbb{P}_2\backslash C\times X\arrow{r}{\tilde{s}_3}[swap]{9} \arrow{d}{X}[swap]{\text{proj}_1} & R_3/M \\
N \arrow{r}{\lambda}[swap]{6} & \mathbb{P}_2\backslash C & R_3 \arrow{u}{\pi}[swap]{M}\arrow{d}{X}[swap]{\Pi}\\
 & \beta^{-1}(0) \arrow{u}{\pi'}[swap]{M}\arrow{r}{\varphi}[swap]{9} &R_3/X\\
Y \arrow{ur}{\lambda'}[swap]{6} \arrow{rr}[swap]{9} \arrow{uu}[swap]{M}&& Y/S\arrow{u}{\Lambda'}[swap]{6}& W' \arrow{uul}{\Phi'}[swap]{6}\arrow{l}{X}
\end{tikzcd}
\]

We shall now prove the following.\\

\e{thm}\label{WW'}
There is a biholomorphism $F:W\to W'$ such that the following diagram commutes.
\[
\begin{tikzcd}
W \arrow{rr}{F}\arrow{dr}[swap]{\Phi}&& W'\arrow{dl}{\Phi'}\\
&R_3&
\end{tikzcd}
\]
\ee{thm}
\e{proof}
The proof mainly consists of unravelling definitions. Recall that
\[
N=\Set{(p,q)\in C\times \mathbb{P}_2\backslash C}{\text{$q$ is on the tangent line to $C$ at $p$}}.
\]
It follows from the definitions of $W$ and $W'$ that 
\begin{align*}
W&=\Set{(f, [(p,q), x])\in R_3\times (N\times X)/S}{\pi(f)=\tilde{s}_3(q,x)},\\
W'&=\Set{(f, [(p, q),h])\in R_3\times (N\times \beta^{-1}(0))/S}{\Pi(f)=\varphi(h), \,q=\pi'(h)}.
\end{align*}
Using Proposition \ref{Cubic Complement}, we can identify $\mathbb{P}_2\backslash C$ with $\chi^{-1}(0)=\beta^{-1}(0)/M$ and define the map
\[
F: W\to W', \quad(f, [(p,q), x])\mapsto (f, [(p, \pi(f\circ \tau_{x})), f\circ\tau_{x}]).
\]
To show that this is a well-defined map into $W'$, we shall firstly show that it is defined independently of the choice of representatives. Let $(f, [(p,q), x])$, $(f, [(p',q'), x'])\in W$ such that $[(p,q), x]=[(p',q'), x']$. Then the action of $S$ on $N\times X$ defined in Section \ref{MThm} is such that for some $t\in S$, $x=x'+t$. Then we have
\begin{align*}
F(f, [(p,q), x])&=(f, [(p, \pi(f\circ \tau_{x})), f\circ\tau_{x}])\\
& =(f, [(p, \pi(f\circ \tau_{x'+t})), f\circ\tau_{x'+t}])\\
& = (f, [(p', \pi(f\circ \tau_{x'})), f\circ\tau_{x'}])\\
& = F(f, [(p',q'), x']).
\end{align*}
Hence $F$ is independent of the representative chosen.

To finish showing that $F$ is well defined, we must show that $F(W)\subset W'$. We shall begin by showing that $f\circ\tau_{x}\in \beta^{-1}(0)$ if $(f, [(p,q), x])\in W$.  Indeed, $\beta(f\circ\tau_{x})=\beta(f)-3x$. Let $\chi: R_3/M\to X$ be the map induced by the $M$-invariant map $\beta$, let $s_3: X\to X$, $x\mapsto 3x$, and let $\tilde{\pi}$ be the projection $W\to N\times X$. Then the following commuting diagram shows that $\beta(f)=3x$.
\[
\begin{tikzcd}
N\times X\arrow{d}{\lambda\times \text{id}_X} \arrow{r}& (N\times X)/S \arrow{d}{\Lambda}&W\arrow{d}{\Phi}\arrow{l}{\tilde{\pi}}\\
\mathbb{P}_2\backslash C\times X \arrow{r}{\tilde{s}_3}\arrow{d}{\proj_2}&R_3/M\arrow{d}{\chi}& R_3\arrow{l}{\pi}\arrow{ld}{\beta}\\
X\arrow{r}{s_3}&X&
\end{tikzcd}
\]
Indeed, the diagram shows that $\beta(f)=\chi(\pi(f))$. Also, since $\pi(f)=\tilde{s}_3(q,x)$, it follows that $\chi(\tilde{s}_3(q,x))=3x$. Hence $f\circ\tau_{x}\in \beta^{-1}(0)$ and $\pi(f\circ \tau_{x})$ can be identified with a point in $\mathbb{P}_2\backslash C$. It follows that $F(W)\subset R_3\times N\times \beta^{-1}(0)$. Checking the conditions defining $W$, obviously $\varphi(f\circ\tau_{x})=\Pi(f)$ and $\pi(f\circ\tau_{x})=\pi(f\circ\tau_{x})$, so $F(W)\subset W'$.

By taking local holomorphic sections of the covering map $N\times X\to (N\times X)/S$ we can show that $F$ is holomorphic by showing that the map
\[
R_3\times N\times X \to R_3\times N\times R_3, \quad(f, (p,q), x)\mapsto (f, (p, \pi(f\circ \tau_{x})), f\circ\tau_{x})
\]
is holomorphic. But this follows from the fact that $X$ acts holomorphically on $R_3$ and the map $\pi$ is holomorphic.

Now, it is clear that $\Phi'\circ F=\Phi$. Since both $W$ and $W'$ are manifolds, to show that $F$ is a biholomorphism, it suffices to show that $F$ is bijective \cite[Proposition 46 A.1]{KnK}. We shall do this by constructing an inverse for $F$ (it is not obvious that this inverse is holomorphic). Since $X$ acts freely on $R_3$, given $(f, (p,q), h)\in R_3\times N\times  \beta^{-1}(0)$ such that $(f,[(p,q),h])\in W'$, there exists a unique $x(f,h)\in X$ such that $h=f\circ\tau_{x(f,h)}$. Furthermore, given $f\in R_3$ and $x\in X$, the construction of the map $\tilde{s}_3$ in the proof of Proposition \ref{Top Cover} shows that $\tilde{s_3}(\pi(f),x)=\pi(f\circ \tau_{-x})$. Then the preimage of $\pi(f)$ under $\tilde{s}_3$ is the set
\[
\Set{(\pi(f\circ\tau_{x_i}), x_i)\in \mathbb{P}_2\backslash C\times X}{i=1,\dots,9},
\]
where $x_1,\dots, x_9$ are the $9$ distinct points on $X$ such that $3x_1=\dots=3x_9=\beta(f)$.

We claim that the map
\[
G: W'\to W, \quad(f, [(p,q'), h])\mapsto (f, [(p, \pi(f\circ\tau_{x(f,h)})), x(f,h)]),
\]
is the inverse of $F$. Indeed, $G$ can be shown to be well defined analogously to $F$. Furthermore,
\begin{align*}
(F\circ G)(f, [(p, q'), h]))&= (f, [(p, \pi(f\circ \tau_{x(f,h)})),f\circ\tau_{x(f,h)}])\\
&= (f, [(p, \pi(h)), h])\\
&= (f, [(p, q'), h]),\\
(G\circ F)(f,[(p,q),x]) &= (f, [(p, \pi(f\circ \tau_{x(f,f\circ\tau_x)})), x(f, f\circ\tau_{x})])\\
&= (f, [(p, q), x]).
\end{align*}
Hence $F$ is a biholomorphism.
\ee{proof}

We shall now explicitly describe the fibres and determine the critical values of $\Phi$ and $\Phi'$. Note that the above theorem implies that it is enough to do the latter for $\Phi'$ only. In order to proceed, we will need to know more about the fibres of $\lambda: N\to \mathbb{P}_2\backslash C$. Let $S^3X$ denote the three-fold symmetric product of $X$ and let
\[
\psi: S^3X\to X, \quad p_1+p_2+p_3\mapsto p_1+p_2+p_3,
\]
be the Jacobi map, where the $+$ sign in $S^3X$ denotes the formal sum of points in $X$, and the $+$ sign in $X$ denotes the Lie group sum in $X$. If $f\in R_3$, then every fibre of $f$, consisting of points $p_1$, $p_2$, $p_3$ counted with multiplicity, has an associated divisor $p_1+p_2+p_3\in S^3X$. When the points $p_1, p_2, p_3$ are not distinct, the fibre is the preimage of a critical value of $f$, and we will call the associated divisor a \emph{branch divisor} of $f$. As above, identify $\mathbb{P}_2\backslash C$ with $\beta^{-1}(0)/M$. Also, since $C$ is the image of an embedding of $X$ in $\mathbb{P}_2$, for each $p\in C$, we shall denote by $x(p)$ the unique point in $X$ corresponding to $p$ under this embedding.\\

\e{prop}\label{Banana}
Let $f\in \beta^{-1}(0)$. Then the branch divisors of $f$ are of the form $x(p)+x(p)+(-2x(p))$, where $(p,q)\in N$ and $q=\pi(f)$.
\ee{prop}

This proposition is an immediate consequence of the following rephrasing of Proposition \ref{LineLemma}.\\

\e{prop}
Let $f\in \beta^{-1}(0)$ and $q=\pi(f)$. Then $D\in \psi^{-1}(0)$ is a divisor corresponding to a fibre of $f$ if and only if $D=x(p_1)+x(p_2)+x(p_3)$, where $p_1, p_2, p_3\in C$ are the intersection points of a line passing through $q$ with $C$, counted with multiplicity.
\ee{prop}

We are now ready to determine the fibres of $\Phi'$. The definition of $W'$ and Proposition \ref{Banana} give the following proposition.\\

\e{prop}\label{bignose}
Let $f\in R_3$, $(p,q)\in N$ and $h\in \beta^{-1}(0)$. Then $(f,[(p,q), h])\in (\Phi')^{-1}(f)$ if and only if all of the following conditions are met.

\e{enumerate}
\item[{\normalfont(1)}] $x(p)+x(p)+(-2x(p))$ is a branch divisor of $h$.
\item[{\normalfont(2)}] $q=\Pi(h)$.
\item[{\normalfont(3)}] There exists $t\in X$ such that $h=f\circ\tau_{-t}$.
\ee{enumerate}
\ee{prop}

The maps $f\in R_3$ whose fibres consist of fewer than $6$ points are those that have a translate $h\in \beta^{-1}(0)$ with fewer than six branch divisors. But then $f$ will have fewer than six branch divisors. By the Riemann-Hurwitz formula, $f$ has a total branching order of $6$. Since $f$ is of degree $3$, this implies that $f$ is a critical value of $\Phi'$ if and only if $f$ has a singleton fibre. By Proposition \ref{Banana}, this corresponds to the case when $q=\pi(h)$ is on the tangent line to an inflection point $p$ of $C$, in which case $(p,q)$ is a branch point of the map $\lambda: N\to \mathbb{P}_2\backslash C$.

If $(f,[(p,q), h])\in (\Phi')^{-1}(f)$, then the equivalence class $[(p,q), h]$ essentially determines a branch divisor of $f$, in that each representative contains the data for a branch divisor of a translate of $f$, and any two representatives in the same class give translates of the same branch divisor. Hence we can loosely interpret the map $\Phi'$ as associating to each $f\in R_3$ its generically six branch divisors.\\

\e{cor}
The critical values of $\Phi$ and $\Phi'$ are those $f\in R_3$ that have a singleton fibre. Equivalently, $f$ has a double critical point.
\ee{cor}

In particular, the set of functions in $R_3$ with a double critical point forms a hypersurface $H$ in $R_3$. This hypersurface, consisting of the critical values of $\Phi$ (or equivalently of $\Phi'$), is invariant under the action of $M$ and the action of $X$. Let $T$ be the set of critical values of $\lambda$. Then $H$ is a `twisted version' of $M\times X\times T$. The hypersurface $H$ is singular because $T$ is singular, being the union of nine distinct lines in $\mathbb{P}_2\backslash C$ with no intersection points on $C$ (see the end of Section \ref{OBC}). Hence we have the following proposition.\\

\e{prop}
The singular points of $H$ are those functions with at least two double critical points.
\ee{prop}

Now, Theorem \ref{WW'} and Proposition \ref{bignose} allow us to explicitly describe the fibres of $\Phi$.\\

\e{prop}
Let $f\in R_3$, $(p,q)\in N$ and $x\in X$. Then $(f, [(p,q), t])\in \Phi^{-1}(f)$ if and only if the following conditions are met.
\e{enumerate}
\item[{\normalfont(1)}] $x(p)+x(p)+(-2x(p))$ is a branch divisor of $f\circ\tau_{t}$.
\item[{\normalfont(2)}] $q=\pi(f\circ\tau_t)$.
\item[{\normalfont(3)}] $3t=\beta(f)$.
\ee{enumerate}
\ee{prop}

The fibres of $\Phi$ can be interpreted in the same way as the fibres of $\Phi'$. Indeed if $(f, [(p,q), t])\in \Phi^{-1}(f)$, then the equivalence class $[(p,q), t]$ corresponds to a branch divisor of $f$, in that every representative provides the data for a translate of this branch divisor.

Finally, it follows from Proposition \ref{Banana} and Section \ref{StrongDom} that $R_3$ contains maps with three double critical points if and only if $X$ is biholomorphic to the equianharmonic cubic.


\renewcommand{\theequation}{\Alph{chapter}.\arabic{section}.\arabic{equation}}

\appendix

\addcontentsline{toc}{chapter}{Bibliography}

\bibliographystyle{plain}

\end{document}